\documentclass[MSNbibl,number,citesort,seceqn,dvips]{arxbj}
\usepackage{dcolumn}
\usepackage[left]{vmkcol}
\usepackage{graphicx}
\aid{0}
\volume{20}
\issue{2}
\pubyear{2014}
\firstpage{645}
\lastpage{675}
\doi{10.3150/12-BEJ501} %kopijuoti is PTS

\makeatletter

\newcommand{\R}{\mathbb{R}}
\newcolumntype{d}[1]{D{.}{.}{#1}}

\newtheorem{lem}[theo]{Lemma}
\newtheorem{theo}{Theorem}[section]
\newremark{remark}{Remark}

\newcommand{\bfSigma}{\boldsymbol{\Sigma}}
\newcommand{\bfB}{\mathbf{B}}
\newcommand{\bfH}{\mathbf{H}}
\newcommand{\Exp}{{\mathbb{E}}}
\newcommand{\RL}{{\mathbb{R}}}
\newcommand{\gr}{\operatorname{gr}}
\newcommand{\argmax}{\operatorname{\arg\max}}
\renewcommand{\emptyset}{\varnothing}
\makeatother

\begin{document}
\begin{frontmatter}

\title{Simple simulation of diffusion bridges with application to
likelihood inference for~diffusions}
\runtitle{Simple simulation of diffusion bridges}

\begin{aug}
%%%% inicialai - be tarpu
\author[1]{\fnms{Mogens} \snm{Bladt}\corref{}\thanksref{1}\ead[label=e1]{bladt@sigma.iimas.unam.mx}} \and
\author[2]{\fnms{Michael} \snm{S{\o}rensen}\thanksref{2}\ead[label=e2]{michael@math.ku.dk}\ead[label=u1,url]{www.math.ku.dk/\textasciitilde michael}}
\runauthor{M. Bladt and M. S{\o}rensen} %% auto
\address[1]{Instituto de Investigacion en Matem\'aticas Aplicadas y en
Sistemas, Universidad Nacional~Aut\'o\-noma de M\'exico, A.P. 20-726,
01000 Mexico, D.F., Mexico. \printead{e1}}
\address[2]{Department of Mathematical Sciences, University of
Copenhagen, Universitetsparken 5, DK-2100 Copenhagen \O, Denmark.
\printead{e2,u1}}
\end{aug}

% HISTORY:
\received{\smonth{4} \syear{2011}}
\revised{\smonth{11} \syear{2012}}

% ABSTRACT
%
\begin{abstract}
With a view to statistical inference for discretely observed diffusion
models, we propose simple methods of simulating diffusion bridges,
approximately and exactly. Diffusion bridge simulation plays a fundamental
role in likelihood and Bayesian inference for diffusion processes.
First a simple method of simulating approximate diffusion bridges is
proposed and studied. Then these approximate bridges are used as
proposal for an easily implemented
Metropolis--Hastings algorithm that produces exact diffusion
bridges. The new method utilizes time-reversibility properties of
one-dimensional diffusions and is applicable to all one-dimensional
diffusion processes with finite speed-measure. One advantage of the
new approach is that simple simulation methods like the Milstein
scheme can be applied to bridge simulation. Another advantage over
previous bridge simulation methods is that the proposed method works
well for diffusion bridges in long intervals because
the computational complexity of the method is linear in the length of
the interval. For $\rho$-mixing diffusions the
approximate method is shown to be particularly accurate for long time
intervals. In a simulation study, we investigate the accuracy and
efficiency of the approximate method and compare it to exact
simulation methods. In the study, our method provides a very good
approximation to the distribution of a diffusion bridge for bridges
that are likely to occur in applications to statistical inference. To
illustrate the usefulness of the new method, we present an
EM-algorithm for a discretely observed diffusion process.
\end{abstract}

% KEYWORDS
% visi is mazosios raides ir pagal abecele
%

\end{frontmatter}

%s1 #&#
\section{Introduction}\label{intro}

In this paper, we propose a simple general method for the simulation of
a one-dimensional diffusion bridge. Our main
motivation is that simulation of diffusion bridges plays
a fundamental role in simulation-based likelihood inference (including
Bayesian inference) for discretely sampled diffusion processes and other
diffusion-type processes like stochastic volatility models.

Our approach is based on the following simple way of constructing a
process that at time zero starts from $a$ and at time $T$ ends in $b$,
where $a$ and $b$ are given numbers. One diffusion process,
$X^{(1)}_t$, is started from the point $a$, while an independent
diffusion, $X^{(2)}_t$, with the same dynamics, is started from the
point $b$. The time of the second diffusion is reversed, so that the
time starts at $T$ and goes downwards to zero. Suppose the process
$X^{(1)}_t$ hits the sample path of the time reversed diffusion
$X^{(2)}_{T-t}$, and let $\tau$ denote the first time the two paths
intersect. Then the process that for $t \leq\tau$ is equal to
$X^{(1)}_t$, and for $t > \tau$ equals $X^{(2)}_{T-t}$, is obviously a
process in the time interval $[0,T]$
that starts at $a$ and ends at $b$. Conditional on the event that the
two processes intersect, we show that the process constructed in this
way is indeed an approximation to a realization of a diffusion bridge
between the two points. A simple rejection sampler is thus obtained by
repeatedly simulating the two diffusions until they hit each other. The
diffusions can be simulated by means of simple procedures like the
Milstein scheme, see Kloeden and Platen \cite{KlPl}, so the new method is easy to
implement for likelihood inference for discretely sampled diffusion
processes. This approximate diffusion bridge is used as proposal for a
Metropolis--Hastings algorithm that has an exact diffusion bridge as
its target distribution. The algorithm uses the pseudo-marginal
approach of Andrieu and Roberts \cite{andrieuroberts} and is easy to implement: to
calculate the rejection probability, a number of independent
diffusions are simulated and it is determined whether or not they
intersect the proposed bridge trajectory.\looseness=1

Diffusion bridge simulation is a highly non-trivial problem that has
been investigated actively over the last 10--15 years. A lucid exposition
of the problems and the state-of-the-art can be found in Papaspiliopoulos and Roberts \cite{OPGR}.
It was previously thought impossible to simulate diffusion bridges by
means of simple procedures, because a rejection
sampler that tries to hit the prescribed end-point for the bridge
(or a small neighbourhood around it) will have an excessively high
rejection probability. The rejection sampler presented in this paper
has a quite acceptable rejection probability because what must be hit
is a sample path rather than a point. The first diffusion bridge simulation
methods in the literature were based on the Metropolis--Hastings algorithm
with a proposal distribution given by a process that is forced to go
from $a$ to $b$, see, for example, Roberts and Stramer \cite{roberts} or Durham and Gallant \cite{durhamgallant}.
Later Beskos, Papaspiliopoulos and Roberts \cite{Beskos2007,beskos7}
developed algorithms for
exact simulation of diffusion bridges. These are rejection sampling
algorithms that use in a clever way a measure change and simulation of
a Brownian bridge, which can easily be simulated. Under strong
boundedness conditions the algorithm is relatively simple, whereas it
is more
complex under weaker condition. Lin, Chen and Mykland \cite{linchenmykland} proposed a
sequential Monte Carlo method for simulating diffusion bridges with a
resampling scheme guided by the empirical distribution of backward
paths. The spirit of this approach has similarities to the methods
proposed here.

An advantage of the method proposed in the present paper is that the same
simple algorithm can be used for all one-dimensional diffusions with a
finite speed-measure, and that it is easy to understand and to
implement. It is also worth noting that the method does not require
that the diffusion is transformed into one with unit diffusion
coefficient. Though such a transformation exists under very general
conditions for one-dimensional diffusions, the transformation is not
in closed form for many interesting examples. Such a transformation
is, for instance, required for the exact algorithm of
Beskos, Papaspiliopoulos and Roberts \cite{Beskos2007,beskos7}. Another,
and perhaps more important,
advantage is that it works particularly well for long time
intervals. For ergodic diffusions, the computational complexity of our
method is shown to be linear in the distance between the two end-points of
the diffusion bridge. This is illustrated in a simulation study where
the computer time increases linearly with the interval length, while
it seems to grow at least exponentially with the interval length for
the exact EA algorithms of Beskos, Papaspiliopoulos and Roberts \cite{Beskos2007}. Thus, the EA algorithm
is likely not to work for long time intervals. Importantly, it is
shown that for ergodic diffusions the approximate method proposed here
simulates an essentially exact
diffusion bridge in long time intervals (apart from the discretization
error). For exponentially mixing diffusions, which covers most
diffusions used in practice, the distribution of the simulated process
goes to that of a diffusion bridge exponentially fast as a function of
the interval length. Thus the proposed method provides a useful
supplement to previously published methods because it works
particularly well for
long time intervals, where the other methods tend not to work. It is
worth noting that simulation-based likelihood inference for discretely
sampled diffusions is mainly important for long time intervals,
because for short time intervals several simpler methods provide
highly efficient estimators, see the following discussion.\looseness=1

The main challenge to likelihood based inference for diffusion models
is that
the transition density, and hence the likelihood function, is not
explicitly available and must therefore be approximated. When the
sampling frequency is relatively high, which is often the case for
financial data, rather crude approximations to the likelihood
functions, like those in Ozaki \cite{ozaki}, Bollerslev and Wooldridge \cite{bollerslev},
Bibby and S{\o}rensen \cite{bmb&ms1} and Kessler \cite{kessler97},
give estimators with a high efficiency. This follows from results
based on high frequency asymptotics in S{\o}rensen \cite{ms10}. When the
interval between the observation times is relatively long, more
accurate approximations to the transition density are needed. One approach
is numerical approximations, either by solving the Kolmogorov PDE
numerically, for example, Poulsen \cite{poulsen} and Hurn, Jeisman and Lindsay \cite{wood&trees}, or by
expansions, for example, A\"{\i}t-Sahalia \cite{ait-sahalia,ait-sahalia2} and
Forman and S{\o}rensen \cite{julie}. Alternatively, likelihood inference can be based on
simulations, an approach that goes back to the seminal paper by
Pedersen \cite{pedersen}, whose method is, however, computationally costly
because he did not use bridge simulation. The inference problem can be
viewed as a missing data problem. If the diffusion
process had been observed continuously, the likelihood functions would
be explicitly given by the Girsanov formula, but the diffusion has
been observed at discrete time points only, so the continuous-time
paths between the observation points can be considered as missing
data. This way of viewing the problem, which goes back to
Dacunha-Castelle and Florens-Zmirou \cite{castellezmirou}, makes it natural to apply either the
EM-algorithm or the Gibbs sampler. To do so, it is necessary to
simulate the missing continuous paths between the observations
conditional on the observations, which is exactly simulation of
diffusion bridges. It was a significant break-through when this was
simultaneously realized by several authors, see Roberts and Stramer \cite{roberts},
Elerian, Chib and Shephard \cite{elerian}, Eraker \cite{eraker} and Durham and Gallant \cite{durhamgallant}, and
approaches based on bridge simulation has since been used by several
authors including Golightly and Wilkinson \cite{golightly&wilkinson1,golightly&wilkinson3,golightly&wilkinson2},
Beskos, Papaspiliopoulos and Roberts \cite{Beskos2006}, Delyon and Hu \cite{delyon&hu}, Beskos, Papaspiliopoulos and Roberts \cite{beskos9} and
Lin \textit{et al}. \cite{linchenmykland}. To illustrate how our bridge simulation
method can be used for likelihood inference, we modify an EM-algorithm
in Beskos, Papaspiliopoulos, Roberts and Fearnhead \cite{Beskos2006} by using our simple simulation
method.\looseness=1

Diffusion bridge simulation is also crucial to simulation-based
inference for other types of diffusion process data than discrete time
observations. Chib, Pitt and Shephard \cite{chibpittshephard} presented a general approach to
simulation-based Bayesian inference for diffusion models when the data
are discrete time observations of rather general, and possibly random,
functionals of the continuous sample path, see also
Golightly and Wilkinson \cite{golightly&wilkinson4}. This approach covers
for instance diffusions observed discretely with measurement error and
discretely sampled stochastic volatility models. In this approach too,
the underlying
continuous time diffusion process must be simulated conditionally on the
observations, which is done by a Metropolis--Hastings algorithm. The
algorithm mixes badly if the entire sample path is updated simultaneously,
so the interval is divided into random subintervals that are updated
sequentially. The sample path in a subinterval must be simulated
conditionally of the values of the diffusion in the other intervals,
which by the Markov property is diffusion bridge simulation given the
values at the end-points. This method can be modified by using the
bridge simulation method proposed here. Baltazar-Larios and S{\o}rensen \cite{fernando} presented an
EM-algorithm for integrated diffusions observed discretely with
measurement error based on the ideas in Chib \textit{et al}. \cite{chibpittshephard}, but
using the bridge simulation method of the present paper, and showed
that the method worked well in simulation studies. We shall briefly
review the results of this paper.\looseness=1

The paper is organized as follows. In Section~\ref{sec2}, we first present the new
approximate bridge simulation method and show in what sense it
approximates a diffusion bridge. Then results are given about exactness
and small rejection probabilities for long time intervals. Finally,
the approximate bridges are used as proposal in a Metropolis--Hastings
algorithm that has an exact diffusion bridge as its target
distribution. In Section~\ref{sim}, the approximate bridge
simulation method is compared to exact simulation methods in two
examples, the Ornstein--Uhlenbeck process and the hyperbolic
diffusion. The study indicates that our method provides a very accurate
approximation to the distribution of a diffusion bridge, except for
bridges that are very unlikely to occur when using the method for
likelihood inference. An EM-algorithm for discretely observed
diffusions based on the proposed bridge simulation method is briefly
presented in Section~\ref{sec4}. It is
demonstrated how the algorithm simplifies for an exponential family of
diffusions (i.e., when drift is linear in the parameters). In this
case, Bayesian inference is considered too. Finally, an application of
the proposed method to estimation for discretely observed integrated
diffusions with measurement errors in Baltazar-Larios and S{\o}rensen \cite{fernando} is briefly
reviewed. Section~\ref{sec5} concludes.

%s2 #&#
\section{Diffusion bridge simulation}\label{sec2}

%s2.1 #&#
\subsection{Approximate bridge simulation}\label{sec2.1}

Let $X = \{ X_t \}_{t\geq0}$ be a one-dimensional diffusion given by the
stochastic differential equation
%
%e2.1 #&#
\begin{equation}
\label{basicsde} \mathrm{d}X_t = \alpha(X_t)\,
\mathrm{d}t + \sigma(X_t)\,\mathrm{d}W_t ,
\end{equation}
where $W$ is a Wiener process, and where the coefficients $\alpha$ and
$\sigma$ are sufficiently regular to ensure that the equation has a
unique weak solution that is a strong Markov process.
Let $a$ and $b$ be given points in the state space of $X$. We present
a method for simulating an approximation to a sample path of $X$ such
that $X_0=a$ and $X_\Delta=b$.
A solution of (\ref{basicsde}) in the interval $[t_1,t_2]$ such that
$X_{t_1}=a$ and $X_{t_2}=b$ will be called a
$(t_1,a,t_2,b)$-bridge. When $t_1=0$ and $t_2=1$, we sometimes
simply call it an $(a,b)$-bridge. We will denote the transition density
of $X$ by $p_t(x,y)$. Specifically, the conditional density of $X_{s+t}$
given $X_s=x$ is $y \mapsto p_t(x,y)$. The state space of $X$ is denoted
by $(\ell,r)$ where $\infty\leq\ell< r \leq\infty$.

Let $W^1$ and $W^2$ be two independent
standard Wiener processes, and define $X^1$ and $X^2$ as the solutions to
\[
\mathrm{d}X^i_t = \alpha\bigl(X^i_t
\bigr)\,\mathrm{d}t + \sigma\bigl(X^i_t\bigr)\,
\mathrm{d}W^i_t , \qquad i=1,2, X^1_0=a
\mbox{ and } X^2_0=b.
\]
The main idea of the paper is to realize an approximation to a
$(0,a,\Delta,b)$-bridge by simulating the process
$X^1$ from $a$ forward in time, and $X^2$ from $b$ backward in
time starting at time $\Delta$. If the samples paths of the two processes
intersect, they can be combined into a realization of a process that
approximates a $(0,a,\Delta,b)$-bridge.

Thus to simulate an approximate diffusion bridge in the interval
$[0,\Delta]$,
we can use any of the several methods available to simulate the
diffusions $X^1$ and $X^2$, see, for example, Kleoden and Platen \cite{KlPl}. Let
$Y^1_{\delta i}$,
$i=0,1, \ldots, N$ and $Y^2_{\delta i}$, $i=0,1, \ldots, N$ be
(independent) simulations of $X^1$ and $X^2$ in $[0,\Delta]$
with step size $\delta= \Delta/N$. Then a simulation of an
approximation to a
$(0,a,\Delta,b)$-bridge is obtained by the following rejection sampling
scheme. Keep simulating $Y^1$ and $Y^2$ until the sample paths cross,
that is, until there is an $i$ such that either $Y^1_{\delta i} \geq
Y^2_{\delta(N-i)}$ and $Y^1_{\delta(i+1)} \leq Y^2_{\delta(N-(i+1))}$
or $Y^1_{\delta i} \leq Y^2_{\delta(N-i)}$ and
$Y^1_{\delta(i+1)} \geq Y^2_{\delta(N-(i+1))}$. Once a
trajectory crossing has been obtained, define
%
%e2.2 #&#
\begin{eqnarray}
\label{bridgesim} B_{\delta i}=
\cases{\displaystyle Y^1_{\delta i} & \quad \mbox{for } $i = 0,1, \ldots, \nu-1$,\vspace*{2pt}
\cr
\displaystyle
Y^2_{\delta(N-i)} & \quad \mbox{for } $i=\nu, \ldots, N,$}
\end{eqnarray}
where $\nu= \min\{ i \in\{1, \ldots, N\} | Y^1_{\delta i} \leq
Y^2_{\delta(N-i)} \}$ if $Y^1_0 \geq Y^2_\Delta$, and
$\nu= \min\{ i \in\{1, \ldots, N\} | Y^1_{\delta i} \geq
Y^2_{\delta(N-i)} \}$ if $Y^1_0 \leq Y^2_\Delta$. Then $B$ approximates
a $(0,a,\Delta,b)$-bridge under the condition of Theorem~\ref{diffusionbridge}. On top of the usual influence of the step size
$\delta$ on the quality of the individual simulated trajectories, the
step size also controls the probability that a trajectory crossing is
not detected. Therefore, it is advisable to choose $\delta$ smaller
than usual.

The rejection probability (the probability of no trajectory crossing)
depends on the drift and diffusion coefficients, on the values of $a$
and $b$, and on the length of the interval $\Delta$. It is shown below
that for ergodic diffusions the rejection probability is small when
$\Delta$ is large (Theorem~\ref{rejectionprob}). Simulation studies in
Section~\ref{sim} indicate that the number of rejections is small when $a$ and
$b$ are not very far apart. When the simulation algorithm is used to
make likelihood
inference for discretely observed diffusion processes (Section~\ref{sec4}), this
is the typical situation for relatively frequent sampling times, and
as just noted, there are in general few rejections when an ergodic diffusion
has been sampled at a low frequency.

The distribution of the process that is simulated by the algorithm
above and the sense in which it is an approximation of a
diffusion bridge is seen from the following theorem, where
%
%e2.3 #&#
\begin{equation}
\label{speed} m(x)=\frac{1}{\sigma^2(x)} \exp \biggl( 2 \int_{z}^x
\frac{\alpha(y)}{\sigma^2(y)} \,\mathrm{d}y \biggr),\qquad  x \in(\ell,r)
\end{equation}
is the density of the speed measure of the diffusion. Here
$z$ is an arbitrary point in the state space $(\ell,r)$.
%
%th2.1 #&#
\begin{theo}
\label{diffusionbridge}
Let $\tau= \inf\{ 0\leq t\leq\Delta| X^1_t=X^2_{\Delta-t} \}$
{ ($\inf\emptyset= +\infty$)} and define
\begin{eqnarray*}
Z_t=\cases{\displaystyle  %
 X^1_t
& \quad \mbox{if } $0\leq t \leq\tau$,\vspace*{2pt}
\cr
\displaystyle X^2_{\Delta-t} & \quad \mbox{if } $\tau< t \leq\Delta$.}
\end{eqnarray*}
Assume that
%
%e2.4 #&#
\begin{equation}
\label{finitespeed} M = \int_\ell^r m(x)\, \mathrm{d}x <
\infty.
\end{equation}
Then the distribution of $\{ Z_t \}_{0\leq t \leq\Delta}$,
conditional on
the event $\{ \tau\leq\Delta\}$, equals the distribution of a
$(0,a,\Delta,b)$-bridge, conditional on the event that the bridge is
hit by an independent diffusion with stochastic differential equation
(\ref{basicsde}) and initial distribution with density $p_\Delta
(b,\cdot)$.
\end{theo}

Before proving Theorem~\ref{diffusionbridge}, we prove a lemma on the
distribution of a time-reversed diffusion. Quite generally, the density
(\ref{speed}) of the speed measure for any one-dimensional diffusion
satisfies the balance equation
%
%e2.5 #&#
\begin{equation}
\label{balance} p_t(x,y)m(x) = p_t(y,x)m(y),
\end{equation}
see Ito and McKean \cite{itomckean}, page 149. Under the condition (\ref{finitespeed})
that the speed measure is finite, an invariant probability measure
exists and has the density function
%
%e2.6 #&#
\begin{equation}
\label{invdensity} \nu(x) = m(x)/M.
\end{equation}

%le2.2 #&#
\begin{lem}\label{lemmafundamental}
Define the time-reversed process $\{ \bar X_t \}$ by $\bar X_t =
X^{2}_{\Delta-t}$. The process $\{ \bar X_t \}$ and the conditional
process $\{ X_t \}$ given that $X_\Delta=b$ have the same transition densities
%
%e2.7 #&#
\begin{equation}
\label{bridgetransition} q(x,s,y,t)= \frac{p_{t-s}(x,y)p_{\Delta-t}(y,b)}{p_{\Delta-s}(x,b)} =
\frac{p_{t-s}(y,x)p_{\Delta-t}(b,y)}{p_{\Delta-s}(b,x)}, \qquad s <
t < \Delta.
\end{equation}
Assume that (\ref{finitespeed}) holds. Then the distribution of $\{
\bar X_t \}$ is equal to the conditional distribution of the process
$\{ X_t \}$ with $X_0 \sim\nu$ given that $X_\Delta=b$.
\end{lem}
\begin{pf}
The second identity in (\ref{bridgetransition}) follows from
(\ref{balance}). The first expression for
$q$ is the well-known expression for the transition density of a
diffusion bridge ending in $b$ at time 1, see Fitzsimmons, Pitman and Yor \cite{yor}, page 111. It
can be easily established by direct calculation.
The second expression for $q$ can similarly be obtained
as the transition density of $\bar X$ by direct calculation. The
conditional density of $\bar X_t$ given $\bar X_s$ ($s<t$) is
\[
p_{\bar X_s, \bar X_t}(x,y)/p_{\bar X_s}(x) = p_{X^2_{\Delta-t},X^2_{\Delta-s}}(y,x)/p_{X^2_{\Delta-s}}(x)
= p_{\Delta-t}(b,y) p_{t-s}(y,x)/p_{\Delta-s}(b,x).
\]

Now suppose that (\ref{finitespeed}) holds, and assume that $X_0 \sim
\nu$. Then $X_\Delta\sim\nu$, and the joint density of $(X_0,
X_\Delta)$ is
$\nu(y_0) p_\Delta(y_0,x) = \nu(x) p_\Delta(x,y_0)$, again by (\ref
{balance}).
Hence, the conditional density of $X_0$ given $X_\Delta=b$ is
$p_\Delta(b,y_0)$.
Obviously, the density of $\bar X_0 = X^2_\Delta$ is $p_\Delta
(b,y_0)$, so the
process $\{ \bar X_t \}$ and the conditional process $\{ X_t \}$
given that $X_\Delta=b$ have the same transition densities and the same
initial distribution. Therefore, they have the same distribution.
\end{pf}
\begin{remark*}
Note that the results of Lemma~\ref{lemmafundamental}
hold for a multivariate diffusion too, provided that a function $v$
exists such
that $p_t(x,y)v(x) = p_t(y,x)v(y)$. Diffusions with this property are
called $v$-symmetric, see the discussion in Kent \cite{kent}. The second
assertion of the lemma holds provided that $v$ is an integrable
function on the state space.
\end{remark*}
\begin{pf*}{Proof of Theorem~\protect\ref{diffusionbridge}}
Let $W^3$ be a
standard Wiener processes independent of $W^1$, and let $X^3$ be the
solution of
\[
\mathrm{d}X^3_t = \alpha\bigl(X^3_t
\bigr)\,\mathrm{d}t + \sigma\bigl(X^3_t\bigr)\,\mathrm{d}W^3_t
,
\]
where the distribution of $X^3_0$ has the density $\nu$ given by
(\ref{invdensity}). Finally, let $\rho$ be the first time
the diffusion $X^3$ hits the sample path of $X^1$. Define a process
by
\begin{eqnarray*}
Y_t=\cases{\displaystyle
X^1_t
& \quad \mbox{if } $0\leq t \leq\rho$,\vspace*{2pt}
\cr
\displaystyle X^3_{t} & \quad \mbox{if } $\rho< t \leq\Delta $}
\end{eqnarray*}
on $\{ \rho\leq\Delta\}$, and $Y = X^1$ on $\{\rho= \infty\}$.
By the strong Markov property $Y$ has the same distribution as $X^1$.
From now on, we condition on $X^3_{\Delta}=b$. Since
\[
P\bigl(Y \in\cdot| X^3_\Delta= b, \rho\leq\Delta\bigr) =
P(Y \in\cdot| Y_\Delta= b, \rho\leq\Delta),
\]
the theorem follows because by Lemma~\ref{lemmafundamental} the
distribution of $\{X^2_{\Delta-t}\}_{0\leq t\leq\Delta}$ equals that of
$\{X^3_{t}\}_{0\leq t\leq\Delta}$ conditional on $X^3_{\Delta}=b$,
so that $P(Y \in\cdot| X^3_\Delta= b, \rho\leq\Delta)
= P(Z \in\cdot| \tau\leq\Delta)$. The event $\{ Y_\Delta= b,
\rho\leq\Delta\}$ is the event that $Y$ is a
$(0,a,\Delta,b)$-diffusion bridge and that the
diffusion bridge is hit by $X^3$, which under the condition
$X^3_{\Delta}=b$ has the initial distribution $p_\Delta(b, \cdot)$
(see the proof of Lemma~\ref{lemmafundamental}).
\end{pf*}

By symmetry, we see that the distribution of the process $\tilde{Z}$
defined by
\begin{eqnarray*}
\tilde{Z}_t=\cases{\displaystyle
X^1_t & \quad \mbox{if } $0\leq t \leq\Delta-\tilde{\tau}$,\vspace*{2pt}
\cr
\displaystyle X^2_{\Delta-t} & \quad \mbox{if } $\Delta- \tilde{\tau} < t \leq
\Delta$}
\end{eqnarray*}
where $\tilde{\tau} = \inf\{ 0\leq t\leq\Delta| X^1_{\Delta
-t}=X^2_{t} \}$,
is that of a $(0,a,\Delta,b)$-bridge conditional on the event that the
bridge is hit by $\{ X^3_{\Delta-t} \}$, where $X^3$ is an independent
diffusion with stochastic differential equation (\ref{basicsde}) and
initial distribution with density $p_\Delta(a,\cdot)$. Here we use $X^1$
until the last time it crosses the trajectory of $\{ X^2_{\Delta-t}\}$,
which happens at time $\Delta-\tilde{\tau}$. Obviously, an approximate
diffusion bridge can also be simulated by using $\tilde{Z}$.

We can consider the diffusions, diffusion bridges and the approximate
diffusion bridge $Z$ as elements of the canonical space, $C_\Delta$,
of continuous functions defined on the time interval $[0,\Delta]$. Each
of these processes induce a probability measure on the usual
sigma-algebra generated by the cylinder sets. Let $f_b$ denote the
Radon--Nikodym derivative of the distribution of the
$(0,a,\Delta,b)$-diffusion bridge with respect to a dominating
measure. The diffusion bridge solves a stochastic differential
equation with the same diffusion coefficient as in (\ref{basicsde}),
see, for example, (4.4) in Papaspiliopoulos and Roberts \cite{OPGR}, so the density $f_b$ is given
by Girsanov's
theorem. Since the drift for the bridge is unbounded at the end point,
one has to choose the dominating measure carefully: it must correspond
to another bridge, see Papaspiliopoulos and Roberts \cite{OPGR}, page 322, and Delyon and Hu \cite{delyon&hu}.
Similarly let $f_a$ and $f_d$ denote the densities of the
distributions of the approximate bridge $Z$ and of a diffusion with
stochastic differential equation (\ref{basicsde}) and initial
distribution with density $p_\Delta(b,x)$, respectively. Let us call a
diffusion of the latter type a $p_\Delta(b)$-diffusion. Finally, for
any $x \in C_\Delta$, let $A_x$ be the set of functions $y \in
C_\Delta$ that intersect $x$. Specifically,
\[
A_x = \bigl\{ y \in C_\Delta| \gr(y) \cap\gr(x) \neq
\emptyset\bigr\},
\]
where $\gr(x) = \{ (t,x_t) | t \in[0,\Delta] \}$. With these
definitions, the relation between the distribution of the approximate
bridge $Z$ and the exact $(0,a,\Delta,b)$-diffusion bridge is
%
%e2.8 #&#
\begin{equation}
\label{z-density} f_a(x) = f_b(x) \pi_\Delta(x)/
\pi_\Delta.
\end{equation}
Here
%
%e2.9 #&#
\begin{equation}
\label{pi} \pi_\Delta(x) = P(Y \in A_x),\qquad
\pi_\Delta= P\bigl((X,Y) \in A\bigr),
\end{equation}
where $A = \{ (x,y) \in C_\Delta^2 | y \in A_x \}$,
$X$ and $Y$ are independent, $X$ is a $(0,a,\Delta,b)$-diffusion
bridge, and $Y$ is a $p_\Delta(b)$-diffusion. Clearly, $\pi_\Delta
(x)$ is the
probability that $Y$ hits the trajectory $x$, while $\pi_\Delta$ is
probability that a $(0,a,\Delta,b)$-bridge is hit by an independent diffusion
with initial distribution $p_\Delta(b, \cdot)$. To prove equation
(\ref{z-density}), note that the joint density of a diffusion bridge
and an independent $p_\Delta(b)$-diffusion given that they intersect
is\looseness=-1
\[
f_b(x) f_d (y) 1_A(x,y)/\pi.
\]\looseness=0
From this expression, (\ref{z-density}) follows by marginalization. Equation
(\ref{z-density}) is important in two ways: it gives an explicit
expression of the quality of our approximate simulation method, and
more importantly, it can be used to improve the approximation.
In\vadjust{\goodbreak}
Section \ref{sec2.3}, we will present two MCMC-algorithms, that improves the
quality of the approximation. In fact, one of them gives exact
diffusion bridges.

Obviously, the quality of our approximate bridge simulation scheme
depends on the probability
$\pi_\Delta$ that a $(0,a,\Delta,b)$-bridge is hit by an independent
diffusion
with initial distribution $p_\Delta(b, \cdot)$. When $\pi_\Delta$
is close
to one, the simulated process is essentially a $(0,a,\Delta
,b)$-bridge. It is
important to realize that the probability $\pi_\Delta$ is not equal
to the
acceptance probability $P(\tau\leq\Delta)$. It is quite possible that
$P(\tau\leq\Delta)$ is small, while $\pi_\Delta$ is close to one. This
happens, for instance, for a diffusion with mean reversion to a level
$\mu$ when $a \ll\mu\ll b$. In the next subsection, we prove that
$\pi_\Delta$ is close to one for long time intervals, provided that the
diffusion is ergodic. In Section~\ref{sim}, we shall investigate when $\pi_\Delta$
can otherwise be expected to be close to one, and when a good approximation
to a diffusion bridge is obtained. Simulations indicate that also
when $\pi_\Delta$ is not close to one (but also not close to zero), the
distribution of the simulated bridge is often indistinguishable from
the distribution of an exact diffusion bridge.

%s2.2 #&#
\subsection{Long time intervals}\label{sec2.2}

We shall now prove that for ergodic diffusions, the probability
$\pi_\Delta$ that a $(0,a,\Delta,b)$-bridge is hit by an independent
diffusion with initial distribution $p_\Delta(b, \cdot)$ is close to
one when $\Delta$ is large. This implies that for large time
intervals, the probability that the process $Z$ defined in Theorem~\ref{diffusionbridge} is a $(0,a,\Delta,b)$-bridge is close to one,
that is, the simulated process is essentially a $(0,a,\Delta,b)$-bridge.
This is very fortunate, because the strength of the method presented
in this paper is that, contrary to other methods for simulating
diffusion bridges, it works well numerically for long intervals, see Section~\ref{sim} and Theorem~\ref{rejectionprob}. A diffusion processes satisfying
(\ref{finitespeed}) is ergodic if
%
%e2.10 #&#
\begin{equation}
\label{infinscale} \int_\ell^{z}
\frac{1}{\sigma^2(x)m(x)} \,\mathrm{d}x = \int_{z}^r
\frac{1}{\sigma^2(x)m(x)} \,\mathrm{d}x = \infty,
\end{equation}
where $m(x)$ is given by (\ref{speed}).

The convergence of the probability $\pi_\Delta$ to one is exponentially
fast in $\Delta$ when the diffusion process has a spectral gap
$\lambda>0$, that is, when the infimum $\lambda$ of the non-zero
eigenvalues of the spectrum of the infinitesimal generator of the
diffusion is strictly positive. For one-dimensional diffusions,
a spectral gap is equivalent to $\rho$-mixing. Most ergodic diffusions
used in practice have this property. Easily checked conditions on the
drift and diffusion coefficients ensuring a spectral gap and hence
exponential convergence of $\pi_\Delta$ can be found in Florens-Zmirou \cite{zmirou},
Hansen and Scheinkman \cite{hansen&scheinkman}, Hansen, Scheinkman and Touzi \cite{hst98} and Genon-Catalot, Jeantheau and Lar\'{e}do \cite{gc00}.

A simple example is a diffusion with linear drift $\alpha(x) =
- \beta(x - \mu)$, $\beta> 0$, which has a spectral gap $\lambda=
\beta$,
provided that the invariant probability measure (\ref{invdensity}) has
finite second moment, cf. Hansen \textit{et al}. \cite{hst98}. Thus by Theorem~\ref{largeintervals}, $1 - \pi_\Delta=
\mathrm{O}(\mathrm{e}^{-\beta\Delta/2})$, so the simple method gives a good
approximation when $\beta\Delta$
is moderately large, that is, when either $\Delta$ is large or when the
diffusion moves fast ($\beta$ large). The Pearson diffusions provide a
useful broad class of diffusions with linear drift; see
Forman and S{\o}rensen \cite{julie}. Ergodic diffusions with linear drift and an arbitrary
stationary distribution were given in Bibby, Skovgaard and S{\o}rensen \cite{bss05}.

For a general diffusion satisfying (\ref{finitespeed}) and
(\ref{infinscale}), Hansen \textit{et al}. \cite{hst98} showed that if the function
\begingroup
\abovedisplayskip=7pt
\belowdisplayskip=7pt
\[
\gamma(x) = \sigma'(x) - \frac{2 \alpha(x)}{\sigma(x)}
\]
has nonzero limits as $x \downarrow\ell$ and $x \uparrow r$, then
the diffusion has a spectral gap $\lambda>0$. Genon-Catalot \textit{et al}. \cite{gc00} gave an explicit
lower bound for $\lambda$, which gives a useful bound on the rate of
convergence of the probability $\pi_\Delta$. To give this lower bound,
we need the functions (defined on the state space of the diffusion)
\[
\varphi(x) = \frac{\int_x^r m(y)\,\mathrm{d}y}{\sigma(x)m(x)}, \qquad \psi(x) = \frac{\int^x_\ell m(y)\,\mathrm{d}y}{\sigma(x)m(x)},
\]
and
\[
C_1(x) = \sup\bigl\{ \varphi^2\bigl(S^{-1}(y)
\bigr) \dvt  y \geq x \bigr\}, \qquad C_0(x) = \sup\bigl\{ \psi^2
\bigl(S^{-1}(y)\bigr) \dvt  y \leq x \bigr\},
\]
where $S$ is the scale function
\[
S(x) = \int_{z}^x \frac{1}{\sigma^2(y) m(y)}\,\mathrm{d}y.
\]
Genon-Catalot \textit{et al}. \cite{gc00} showed that
\[
\lambda\geq C = \frac{1}{8 \inf_{x \in\R} \max\{ C_1(x),C_0(x) \}}.
\]
If the limits of $\gamma(x)$ exist and are nonzero, then $C>0$, and
Theorem~\ref{largeintervals} implies that $1 - \pi_\Delta
= \mathrm{O}(\mathrm{e}^{- C \Delta/2})$.

The results discussed here are summarized in the following theorem.
%
%th2.3 #&#
\begin{theo}\label{largeintervals}
For ergodic diffusions
\[
\pi_\Delta\rightarrow1
\]
as $\Delta\rightarrow\infty$. If the diffusion has a spectral gap
$\lambda>0$, then
\[
1 - \pi_\Delta= \mathrm{O}\bigl(\mathrm{e}^{-\lambda\Delta/2}\bigr).
\]
\end{theo}\endgroup

\begin{pf}
$\pi_\Delta$ is the probability that a $(0,a,\Delta,b)$-bridge
$X^b$ is hit by an independent diffusion $X^d$ with initial
distribution $p_\Delta(b, \cdot)$. Let $x$ denote the point from which
$X^d$ starts, and assume that $x > a$. The case $x < a$ can be treated
similarly.

By Karlin and McGregor \cite{karlinmcgregor} (page 1144) the probability that $X^d_{\Delta
/2} >
q$ and $X^b_{\Delta/2} \leq q$ without having been coincident in
$[0,\Delta/2]$
(conditional on $X^d_0 = x$) is
\[
P_{q,\Delta} = \det\left\{ %
 {P
\bigl(X^b_{\Delta/2} \leq q\bigr) \atop P\bigl(X^d_{\Delta/2} \leq q | X^d_0 =
x\bigr)}
\enskip
{P\bigl(X^b_{\Delta/2}
> q\bigr) \atop P\bigl(X^d_{\Delta/2} > q | X^d_0
= x\bigr) } \right\}.\vadjust{\goodbreak}
\]
Here $q$ is an arbitrary real number. That the result holds for
time-inhomogeneous diffusions too follows from Karlin \cite{karlin}. Since
the diffusion is ergodic, $X^d_{\Delta/2}$ converges weakly to the
stationary distribution with density $\nu$. Then the density
of $X^b_{\Delta/2}$ (cf. (\ref{bridgetransition})) satisfies
\[
q(a,0,y,\Delta/2)= \frac{p_{\Delta/2}(a,y)p_{\Delta/2}(y,b)}{p_{\Delta}(a,b)} \rightarrow\frac{\nu(y)\nu(b)}{\nu(b)} = \nu(y),
\]
as $\Delta\rightarrow\infty$. Hence
\[
P_{q,\Delta} \rightarrow\det\left\{ %
{P_\nu\bigl((-\infty,q]\bigr) \atop  P_\nu\bigl((-\infty,q]\bigr)}
\enskip
{P_\nu\bigl((q,\infty]\bigr) \atop     P_\nu\bigl((q,\infty]\bigr) }
\right\} = 0,
\]
as $\Delta\rightarrow\infty$, where $P_\nu$ denotes the stationary
distribution. This implies that $\pi_\Delta\rightarrow1$.

Now assume that $\lambda> 0$, where $\lambda$ is the
infimum of the nonzero eigenvalues of the spectrum of the
infinitesimal generator of the diffusion. Then $p_{\Delta/2}(x,y) =
\nu(y)(1+\mathrm{O}(\mathrm{e}^{-\lambda\Delta/2}))$, see, for example, Karlin and Taylor
\cite{karlintaylor},
page 332, and it follows that the density function of
$X^b_{\Delta/2}$ also satisfies that $q(a,0,y,\Delta/2) =
\nu(y)(1+\mathrm{O}(\mathrm{e}^{-\lambda\Delta/2}))$. Hence,
\[
P_{q,\Delta} = \det\left\{ %
{
P_\nu\bigl((-\infty,q]\bigr) \atop  P_\nu\bigl((-\infty,q]\bigr)}
\enskip
{P_\nu\bigl((q,\infty]\bigr) \atop P_\nu\bigl((q,\infty]\bigr) }
\right\} + \mathrm{O}\bigl(\mathrm{e}^{-\lambda\Delta/2}\bigr) = \mathrm{O}
\bigl(\mathrm{e}^{-\lambda\Delta/2}\bigr).
\]\upqed
\end{pf}

If we replace $\Delta/2$ by the time $\gamma\Delta$ ($0 < \gamma<
1$) in the
proof, it follows that $1 - \pi_\Delta= \mathrm{O}(\mathrm{e}^{-\gamma\lambda\Delta
})$. In
practice, this refinement does not make much difference.

A similar result can be proved for the rejection probability in the
same way. Let $p_\Delta= P(\tau> \Delta)$ denote the rejection probability.
%
%th2.4 #&#
\begin{theo}\label{rejectionprob}
For ergodic diffusions
\[
p_\Delta\rightarrow0
\]
as $\Delta\rightarrow\infty$. If the diffusion has a spectral gap
$\lambda>0$, then
\[
p_\Delta= \mathrm{O}\bigl(\mathrm{e}^{-\lambda\Delta/2}\bigr).
\]
\end{theo}

Theorem~\ref{rejectionprob} implies that for ergodic diffusions the
computational complexity of our method is linear in the time distance
between the two end-points of the diffusion bridge, provided that the
diffusion is simulated by a scheme that is linear in the interval
length. This is true of simple simulation method like the Milstein
scheme.

%s2.3 #&#
\subsection{Exact bridge simulation}\label{sec2.3}

Since we know the relationship between the distribution of an exact
diffusion bridge and the distribution of our simple approximation,
cf. (\ref{z-density}), it is natural to consider a Metropolis--Hastings
algorithm for which the proposal is our simple simulation method and
the target distribution is the distribution of an exact diffusion
bridge. In each step, we draw an independent sample path $X^{(i)}$ with
distribution given by $f_a$ by means of the simple simulation
method. The proposed sample path is accepted with probability
$\alpha(X^{(i-1)}, X^{(i)}) = \min(1, r(X^{(i-1)}, X^{(i)}) )$,
where
\[
r\bigl(X^{(i-1)}, X^{(i)}\bigr) = \frac{f_b(X^{(i)})f_a(X^{(i-1)})}{f_b(X^{(i-1)})f_a(X^{(i)})} =
\frac{\pi_\Delta(X^{(i-1)})}{\pi_\Delta(X^{(i)})}.
\]
Here $X^{(i-1)}$ is the previously accepted sample path, and
$\pi_\Delta(x)$ is the probability that the sample path $x$ is hit by an
independent $p_\Delta(b)$-diffusion, given by (\ref{pi}). The MH-algorithm
produces draws of exact diffusion bridges, but as $\pi_\Delta(x)$ is
not explicitly known, the algorithm cannot be used as it stands.

One possibility is to use a Monte Carlo within Metropolis
algorithm. This can be done by simulating in each step $N$ independent
$p_\Delta(b)$-diffusions, $ \mathbf{Y}^{(i)} = (Y^{(i,1)}, \ldots,
Y^{(i,N)})$ and
estimate $\pi_\Delta(x)$ consistently by
\[
\tilde\pi_\Delta\bigl(x;{\mathbf{Y}}^{(i)}\bigr) = \frac1N
\sum_{j=1}^N 1_{A_x}
\bigl(Y^{(i,j)}\bigr).
\]
If $N$ is sufficiently large, this would certainly produce a very good
approximation to a diffusion bridge, but not exact diffusion
bridges. Usually, the density function $p_\Delta(b,x)$ is not
explicitly known, but a $p_\Delta(b)$-diffusion can be simulated as
follows. Simulate the solution $V$ to (\ref{basicsde}) in the time
interval $[0,2\Delta]$ with $V_0=b$. Then $Y_t = V_{t+\Delta}$, $t
\in
[0,\Delta]$, is a $p_\Delta(b)$-diffusion.

In order to simulate exact diffusion bridges, we propose an MCMC
algorithm of the pseudo-marginal type studied in Andrieu and Roberts \cite{andrieuroberts}.
The basic idea of the pseudo-marginal approach is to replace the
factor in the acceptance ratio, which we cannot calculate,
$f_b(x)/f_a(x) = 1/\pi_\Delta(x)$ by an unbiased MCMC
estimate. The beauty of the method is that by including the MCMC draws
needed for the estimate of $1/\pi_\Delta(x)$ in the MH-Markov chain,
the marginal equilibrium distribution of the bridge draws is exactly
$f_b$, irrespective of the randomness of the estimate of
$1/\pi_\Delta(x)$.

Define a
random variable $T$ in the following way. For a given sample path $x
\in C_\Delta$, simulate a sequence of independent $p_\Delta(b)$-diffusions
$Y^{(1)},Y^{(2)},\ldots$ until a sample path is obtained that intersects
$x$. Then $T$ is defined as the number of the first $Y^{(i)}$ that hits $x$:
\[
T = \min\bigl\{i\dvt  Y^{(i)} \in A_x \bigr\}.
\]
By results for the geometric distribution $E(T) = 1/\pi_\Delta(x)$,
so if ${\mathbf{T}} = (T_1,\ldots,T_N)$ is a vector of $N$ independent
draws of $T$, then an unbiased and consistent estimator of
$1/\pi_\Delta(x)$ is\looseness=-1
\[
\hat\rho_\Delta(x;{\mathbf{T}}) = \frac1N \sum
_{j=1}^N T_j.
\]\looseness=0
Now consider the following MH-algorithm, where draws of ${\mathbf{T}}$
are included in the Markov chain. In each step, we draw (independently
of the previous draws) a sample paths $X^{(i)}$ by the
simple algorithm, and with $x=X^{(i)}$ we draw a vector of $N$ independent
$T$-values, ${\mathbf{T}}^{(i)} = (T_1^{(i)},\ldots,T_N^{(i)})$. Note that
the distribution of ${\mathbf{T}}^{(i)}$ depends on $X^{(i)}$. The proposed
update $(X^{(i)},{\mathbf{T}}^{(i)})$ is accepted with probability
\[
\hat\alpha\bigl(X^{(i-1)}, {\mathbf{T}}^{(i-1)}, X^{(i)},
{\mathbf{T}}^{(i)}\bigr) = \min\bigl(1, \hat r\bigl(X^{(i-1)},{
\mathbf{T}}^{(i-1)}, X^{(i)}, {\mathbf {T}}^{(i)}\bigr)
\bigr),
\]
where
\[
\hat r\bigl(X^{(i-1)}, {\mathbf{T}}^{(i-1)}, X^{(i)}, {
\mathbf{T}}^{(i)}\bigr) = \frac{\hat\rho_\Delta(X^{(i)};{\mathbf{T}}^{(i)})} {
\hat\rho_\Delta(X^{(i-1)};{\mathbf{T}}^{(i-1)})}.
\]
By results in Andrieu and Roberts \cite{andrieuroberts}, the target distribution of $X$ is
that of an exact diffusion bridge. In fact, since
\[
\hat r \bigl(x^{(1)},{\mathbf{t}}^{(1)},x^{(2)},{
\mathbf{t}}^{(2)}\bigr) = \frac{f_a(x^{(2)}) {\mathbf{f}}_g({\mathbf{t}}^{(2)} | x^{(2)})
\hat\rho_\Delta(x^{(2)};{\mathbf{t}}^{(2)})
f_a(x^{(1)}) {\mathbf{f}}_g({\mathbf{t}}^{(1)} | x^{(1)})} {
f_a(x^{(1)}) {\mathbf{f}}_g({\mathbf{t}}^{(1)} | x^{(1)})
\hat\rho_\Delta(x^{(1)};{\mathbf{t}}^{(1)})
f_a(x^{(2)}) {\mathbf{f}}_g({\mathbf{y}}^{(2)} | x^{(2)})},
\]
where ${\mathbf{f}}_g({\mathbf{t}} | x)$ is the conditional density of
${\mathbf{T}}$ given $X=x$, we see that the density of the target
distribution is
\[
p(x,{\mathbf{t}}) = f_a(x) {\mathbf{f}}_g({\mathbf{t}}
| x) \hat \rho_\Delta (x,{\mathbf{t}}) \pi_\Delta =
f_b(x) {\mathbf{f}}_g({\mathbf{t}} | x) \hat
\rho_\Delta (x,{\mathbf{t}})\pi_\Delta(x),
\]
where we have used (\ref{z-density}). Since $\hat\rho_\Delta
(x,{\mathbf{t}})$ is an unbiased estimator of $1/\pi_\Delta(x)$
conditionally on $x$, we find by marginalizing that the density of $X$
is $f_b$, the density function of an exact diffusion bridge.

To produce diffusion bridges by the proposed pseudo-marginal MH-algorithm,
a number of sample paths of ordinary diffusions must be simulated. If
these sample paths are simulated by an approximate method, like the Milstein
scheme, a small discretization error is introduced. This problem can,
however, be avoided by using the methods for exactly simulating diffusions
developed by Beskos, Papaspiliopoulos and Roberts \cite{Beskos2007,beskos7}.
By combining our exact MH bridge simulation algorithm with exact
diffusion simulation methods, exact diffusion bridges can be
efficiently simulated even in long time intervals.

For ergodic diffusions, the computational complexity of the exact
algorithm is linear in the interval length $\Delta$. In the previous
subsection, we saw that simulation of the proposal is linear in
$\Delta$, and the expected number of $p_\Delta(b)$-diffusions
simulated in the $i$th M--H step is $N/\pi_\Delta(X^{(i)})$. For a
diffusion bridge $X$, the expectation of $\pi_\Delta(X)$ is
$\pi_\Delta$, which by Theorem~\ref{largeintervals} converges to one
as $\Delta\rightarrow\infty$. Since $\pi_\Delta$ is usually not
small, each iteration of the M--H-algorithm producing exact diffusion
bridges is not expected to require much more computing time than the
approximate algorithm. Finally, since the acceptance ratio tends to
one as $\Delta\rightarrow\infty$, the acceptance rate is high
when $\Delta$ is large.

%s3 #&#
\section{Simulation study}
\label{sim}

In this section, we investigate simulation of two examples of
diffusion bridges for which our new methods can be compared to other
exact algorithms, the Ornstein--Uhlenbeck process
and the hyperbolic diffusion. We compare the distribution of exact
bridge simulations to the distribution of the process\vadjust{\goodbreak} obtained by our
approximate bridge simulation method. It is
found that the approximation is very accurate in most cases, including
bridges that are likely to occur in applications to likelihood
inference. We also compare CPU execution times, and illustrate how our
exact M--H simulation method provides an exact bridge in the extreme cases
where the approximate method does not provide a good approximation.

%s3.1 #&#
\subsection{The Ornstein--Uhlenbeck bridge}\label{sec3.1}

First, we consider an Ornstein--Uhlenbeck bridge, which is
a solution to the stochastic differential equation
\[
\mathrm{d}X_t = -\theta X_t \,\mathrm{d}t + \sigma \,\mathrm{d}W_t
\]
conditionally on $X_0=a$ and $X_1=b$ for some $a,b\in\RL$.
From the well-known Gaussian transition densities of the
Ornstein--Uhlenbeck process we can calculate the transition
densities of the Ornstein--Uhlenbeck bridge by (\ref{bridgetransition}).
Thus, we could in principle simulate the Ornstein--Uhlenbeck bridge by
sampling transitions from these densities. The following well-known
alternative method is, however, numerically more stable.
%
%le3.1 #&#
\begin{lem}\label{lemmamichael}
Generate $X_{t_0},X_{t_1}, \ldots, X_{t_n},X_{t_{n+1}}$,
where $0=t_0 < t_1 < \cdots< t_n < t_{n+1}$, by $X_0 = x_0$ and
\[
X_{t_i} = \mathrm{e}^{-\theta(t_i - t_{i-1})}X_{t_{i-1}} + W_i, \qquad i =1,
\ldots, n+1,
\]
where the $W_i$s are independent and $W_i \sim N ( 0,
\sigma^2 (1 - \mathrm{e}^{-2 \theta(t_i - t_{i-1})}  )/(2 \theta)
)$.
Define
\[
Z_{t_i} = X_{t_i} + (x - X_{t_{n+1}})\frac{\mathrm{e}^{\theta t_i}-\mathrm{e}^{-\theta t_i}} {
\mathrm{e}^{\theta t_{n+1}}-\mathrm{e}^{-\theta t_{n+1}}},\qquad
i=0, \ldots, n+1.
\]
Then $(Z_{t_0}, Z_{t_1}, \ldots, Z_{t_n}, Z_{t_{n+1}})$ is
distributed like an Ornstein--Uhlenbeck bridge with $Z_{t_0}=x_0$ and
$Z_{t_{n+1} = x}$.
\end{lem}
%
%f1 #&#
\begin{figure}

\includegraphics{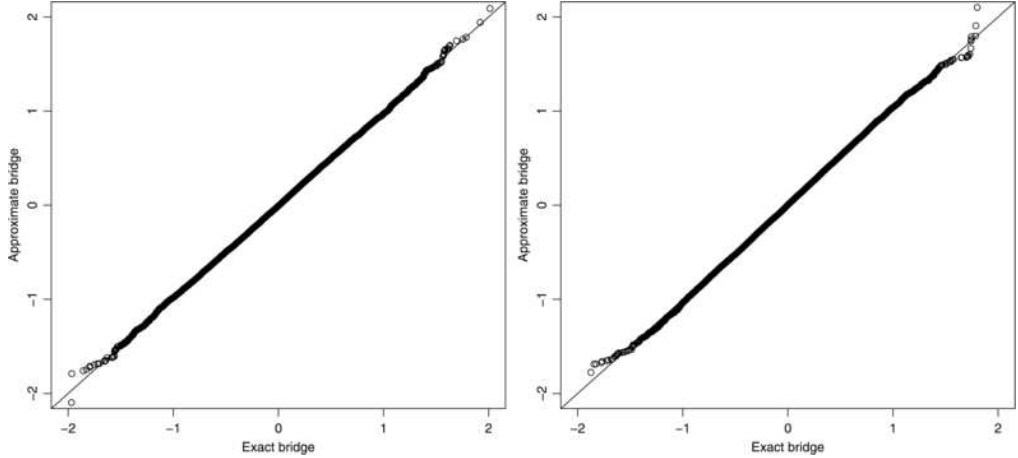}

\caption{Q--Q plots that compare the empirical distribution at time 0.5
based on 25,000 simulated $(0,0)$ diffusion bridges obtained
by our approximate method to
that based on 25,000 exactly simulated diffusion bridges. The left
plot is for the Ornstein--Uhlenbeck bridge and the right plot is for the
hyperbolic diffusion bridge. Exact simulations are obtained by the
method in Lemma \protect\ref{lemmamichael} for the Ornstein--Uhlenbeck bridge
and by the exact algorithm of Beskos, Papaspiliopoulos and Roberts for
the hyperbolic diffusion bridge.}\label{simulation1}
\end{figure}
%
%f2 #&#
\begin{figure}

\includegraphics{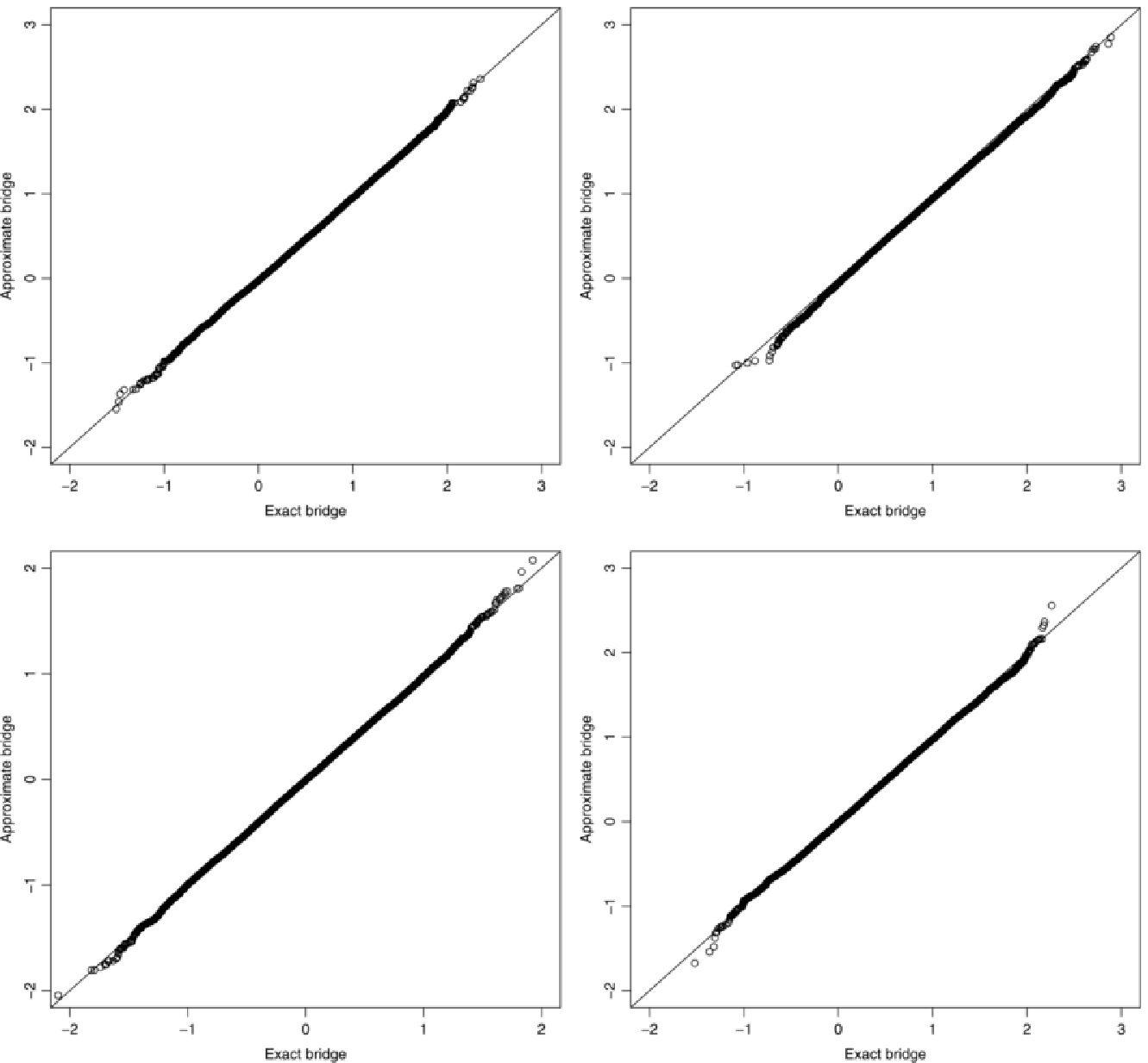}

\caption{Q--Q plots that compare the empirical distributions at time 0.5
based on 25,000 simulated $(0,1)$, $(0,2)$, $(-1,1)$ and $(-1,2)$
Ornstein--Uhlenbeck bridges obtained by our approximate method to that
based on 25,000 exactly simulated Ornstein--Uhlenbeck bridges. Exact
simulations are obtained by the method in Lemma \protect\ref{lemmamichael}.}\label{simulation2}
\end{figure}

In all examples considered in the following we simulated 25,000
(or 10,000) realizations of diffusion bridges over the time interval
$[0,1]$. The Euler scheme was used with discretization level $N=100$ (step
size $\delta=0.01$). For the Ornstein--Uhlenbeck process the Euler
scheme is equal to the Milstein scheme. The methods were implemented
in Fortran 90 on a Dell Precision M65 workstation (laptop).

For the Ornstein--Uhlenbeck bridge, we chose the parameter values
$\theta= 0.5$ and $\sigma=1.0$. First we considered a bridge that
started at $0$ and ended at $0$. We compare our approximate method
based on Theorem~\ref{diffusionbridge} to the exact
algorithm of Lemma~\ref{lemmamichael}. To the left in Figure~\ref{simulation1}, we have plotted the quantiles of the empirical
distribution at the time point $0.5$ obtained by our approximate method
against the quantiles of the empirical distribution obtained by the exact
algorithm. The two distributions appear to be equal.
Similar comparisons of quantiles at time 0.5 for our approximate
method to quantiles of an exact bridge are
presented in Figure~\ref{simulation2} for $(0,1)$, $(0,2)$, $(-1,1)$
and $(-1,2)$ Ornstein--Uhlenbeck bridges. In all four cases, the two
distributions seem to be essentially equal, except for a very small
negative bias for the $(0,2)$-bridge. Similar results were found for
several other comparisons of distributions with similar values of the
start and end points, $a$ and $b$.

The CPU execution time (in seconds) to simulate 10,000 Ornstein--Uhlenbeck
bridges using our approximate method for the various starting points,
$a$, and end points, $b$, are given in Table~\ref{OUtab} together with
estimated rejection probabilities ($p_\Delta$, see Theorem~\ref{rejectionprob}). The table also gives the probabilities that an
Ornstein--Uhlenbeck process moves from $a$ to $b$ or farther in the
time interval $[0,1]$. We see that for moves that are likely to
appear in data sets, the CPU times and rejection probabilities are
small, and the CPU times are only slightly larger than the execution
time for the exact algorithm which is about 0.5 CPU seconds. For more
unlikely moves the rejection probability is quite large, but also in
these cases the execution time is not a problem in applications. The
last column of Table~\ref{OUtab} gives the (estimated) probability of
the event that an exact $(a,b)$-bridge is not hit by an independent diffusion
with initial distribution $p_1(b, \cdot)$. These probabilities were
found by simulating exact Ornstein--Uhlenbeck bridges and independent
Ornstein--Uhlenbeck processes with initial distribution $p_1(b, \cdot)$.
If this probability were zero, our approximate method would simulate
an exact diffusion bridge. The probabilities are small, but not
negligible. It is remarkable that the approximate method gives a quite
accurate approximation to a diffusion bridge in spite
of this. The reason must be that the diffusion bridges are not hit by
the independent diffusion in a systematic way for the $a$ and $b$ values
considered here.
%
%t1 #&#
\begin{table}
\tablewidth=\textwidth
\tabcolsep=0pt
\caption{The CPU execution time (in seconds) used to simulate 10,000
Ornstein--Uhlenbeck bridges using our approximate method for various
starting points, $a$, and end points, $b$. Estimated rejection
probabilities ($p_\Delta$, see Theorem \protect\ref{rejectionprob})
and the probabilities of a move from $a$ to $b$ or farther are
listed too. The last column gives the probability that an exact
$(a,b)$-bridge is not hit by an independent diffusion
with initial distribution $p_1(b, \cdot)$}
\label{OUtab}
\begin{tabular*}{\textwidth}{@{\extracolsep{\fill}}ld{2.1}d{1.2}d{1.4}d{1.2}@{}}\hline
\multicolumn{1}{l}{$a \mapsto b$} & \multicolumn{1}{l}{CPU (sec.)} & \multicolumn{1}{l}{Rejection
prob.}
& \multicolumn{1}{@{}l}{Probability of move} & \multicolumn{1}{l@{}}{$ 1 - \pi$} \\ \hline
$\hphantom{-}{0 \mapsto0}$ & 0.5 & 0.17 & & 0.28 \\
$ \hphantom{-}{0 \mapsto1}$ & 0.7 & 0.41 & 0.1 & 0.21 \\
$ \hphantom{-}{0 \mapsto2}$ & 1.7 & 0.77 & 0.006 & 0.08 \\
$-1 \mapsto1$ & 1.9 & 0.80 & 0.02 & 0.16 \\
$-1 \mapsto2$ & 11.9 & 0.97 & 0.0005 & 0.06 \\ \hline
\end{tabular*}
\end{table}

\begin{figure}[b]

\includegraphics{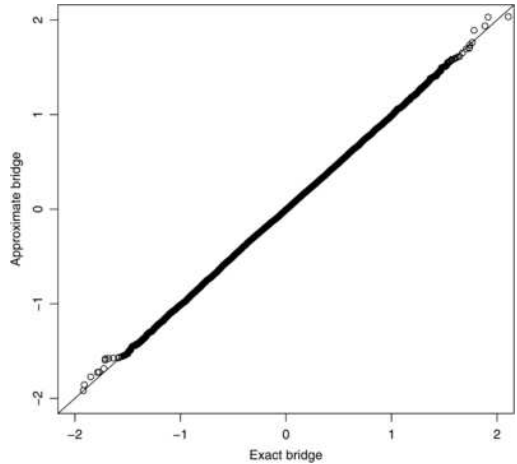}

\caption{Q--Q plot that compares the empirical distribution at time 0.5
based on 25,000 simulated $(-2,2)$ Ornstein--Uhlenbeck bridges obtained
by our approximate method to that based on 25,000 exactly simulated
Ornstein--Uhlenbeck bridges.}\label{simulationextreme}
\end{figure}

In order to test our approximate method in an extreme situation,
we simulated 25,000 Ornstein--Uhlenbeck bridges that started from
$-$2 and ended in 2. The probability that an Ornstein--Uhlenbeck process
with parameters $\theta= 0.5$ and $\sigma=1.0$ moves from $-$2 to 2 or
farther in the time interval $[0,1]$ equals the probability
that a standard normal distribution is
larger than 4.04, which is 0.00003, so this a indeed
a very extreme event. Not surprisingly that rejection
rate was very high, but as appears from Figure~\ref{simulationextreme}
the distribution at time $0.5$ fits the distribution
obtained by exact simulation very well.
%
%f3 #&#

%
%f4 #&#
\begin{figure}

\includegraphics{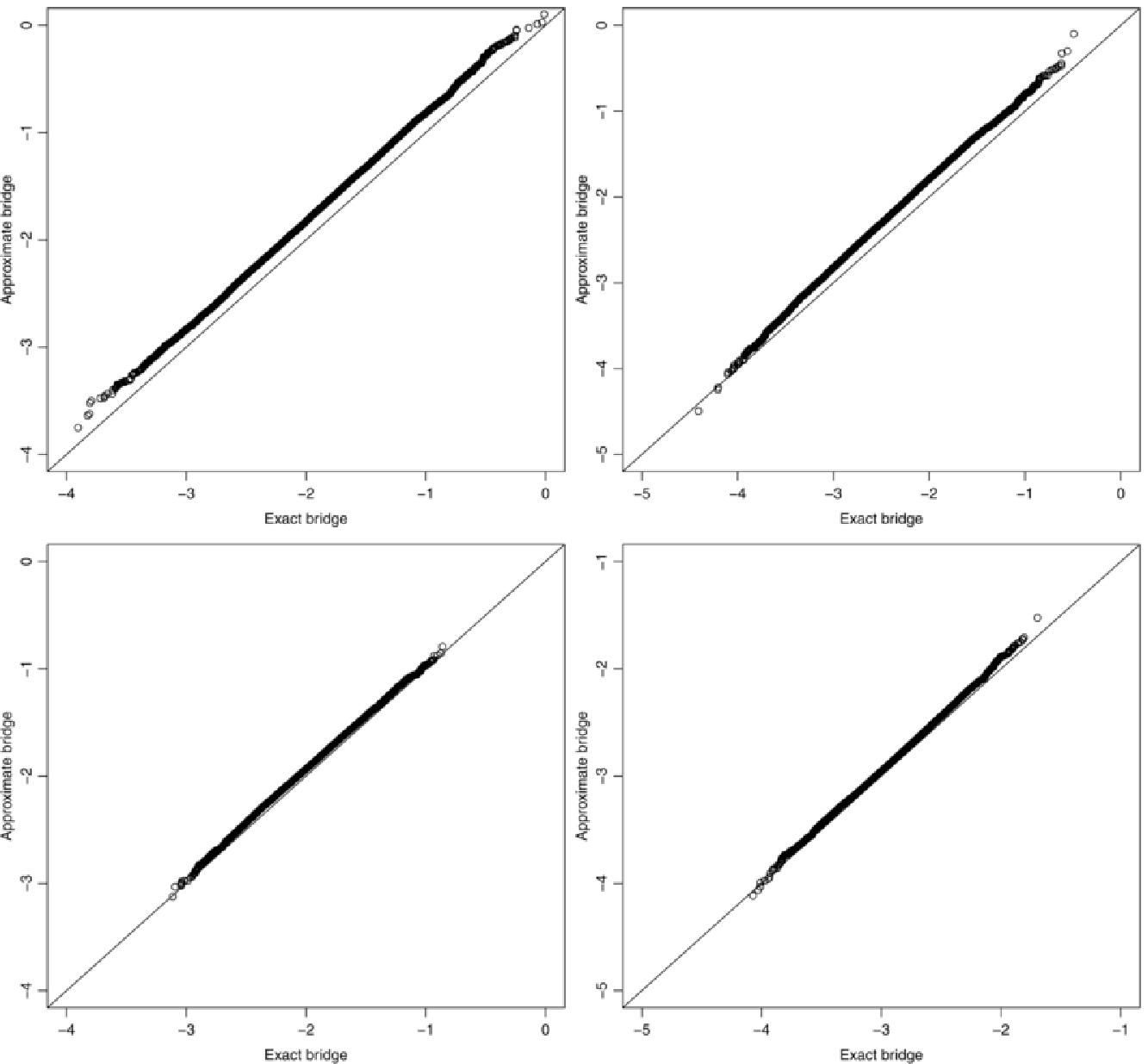}

\caption{Q--Q plots that compare the empirical distributions at time 0.5 (first
row) and time 0.1 (second row) based on 25,000 simulated $(-2,-2)$
(left) and $(-3,-2)$ (right) Ornstein--Uhlenbeck bridges obtained
by our approximate method to that based on 25,000 exactly simulated
Ornstein--Uhlenbeck bridges. Exact simulations are obtained by the
method in Lemma \protect\ref{lemmamichael}.}\label{simulation3}
\end{figure}

%
%f5 #&#
\begin{figure}

\includegraphics{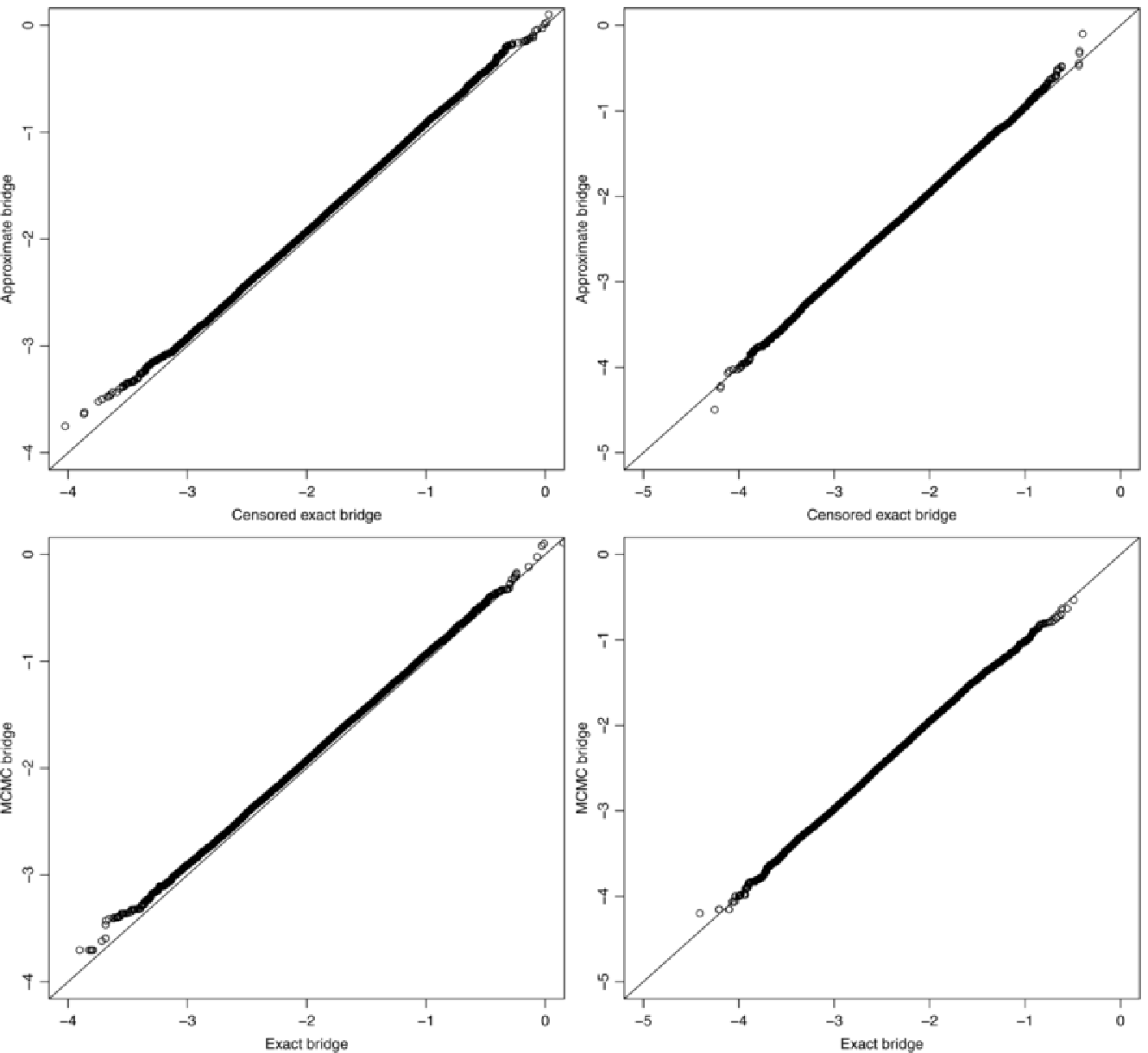}

\caption{In the first row, Q--Q plots compare the empirical distributions at time 0.5
based on 25,000 simulated $(-2,-2)$ (left) and $(-3,-2)$ (right)
Ornstein--Uhlenbeck bridges obtained by our approximate method to
that based on 25,000 exactly simulated Ornstein--Uhlenbeck bridges,
where the bridges that were not hit by an independent diffusion with initial
distribution $p_1(b, \cdot)$ were removed from the sample. In the
second row, Q--Q plots compare the empirical distributions at time 0.5
based on 25,000 simulated $(-2,-2)$ (left) and $(-3,-2)$ (right)
Ornstein--Uhlenbeck bridges obtained by our Metropolis--Hastings
algorithm (after a burn-in of 5000 iterations) to
that based on 25,000 exactly simulated Ornstein--Uhlenbeck bridges.}\label{simulation4}
\end{figure}

The only situation we have been able to find where the distribution
obtained by our approximate simulation method differs appreciably from the
distribution of an exact bridge is when the start and end points, $a$
and $b$, have the same sign and are both far from the equilibrium point
zero. This is to be expected because we simulate an exact bridge
conditional on the event that it is hit by an independent diffusion with
initial distribution $p_1(b, \cdot)$. When $b$ is far from zero, most
of the probability mass of $p_1(b, \cdot)$ is located considerably
closer to zero than $b$ (because of the drift towards zero). The
independent diffusion will tend to move towards zero, while the
$(a,b)$-bridge will tend to stay relatively close to $a$ and
$b$. Only trajectories of the bridge that move sufficiently towards
zero has a reasonable chance of being hit by the independent
diffusion. This creates a bias towards zero. The comparison of
quantiles at time 0.5 and time 0.1 for $(-2,-2)$ and $(-3,-2)$
bridges are presented in Figure~\ref{simulation3}. As expected from the
consideration above, there is a positive bias.
The time point 0.5 was chosen in most simulations because it is
expected that this is the time point where it is most difficult to
get a good approximation of the distribution of a diffusion
bridge. The comparison of distributions at time 0.1 illustrate that
the approximate method works better close to the end-points of
the time interval, where the bias is considerably smaller. Several similar
comparisons confirm that the approximate method works better close to the
end-points than at time 0.5.

%
%f6 #&#
\begin{figure}

\includegraphics{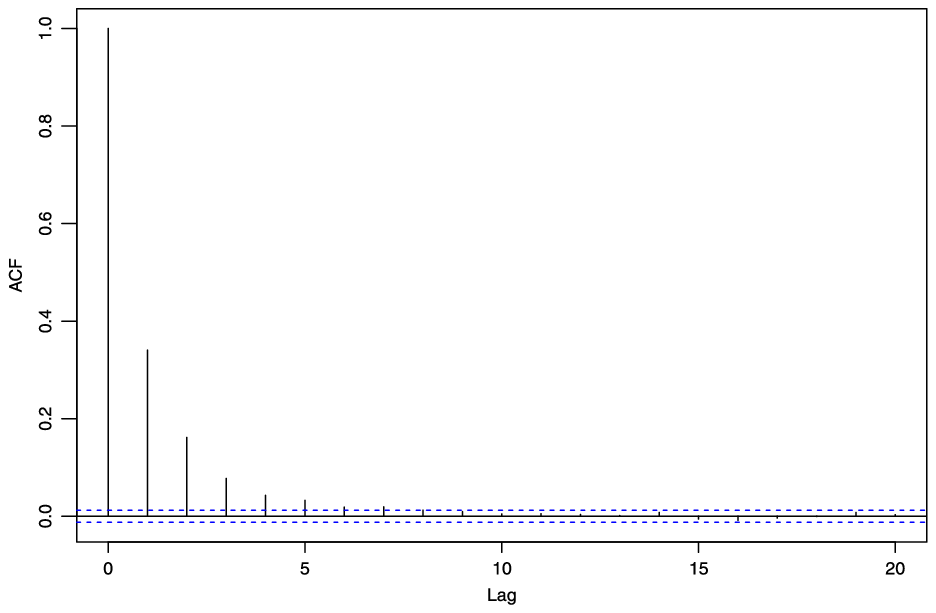}

\caption{The empirical autocorrelation function of the successive values
of the $(-3,-2)$ Ornstein--Uhlenbeck bridge at time 0.5 obtained by
our Metropolis--Hastings algorithm.}\label{simulation5}
\end{figure}

%
%t2 #&#
\begin{table}[b]
\tablewidth=\textwidth
\tabcolsep=0pt
\caption{Estimated rejection probabilities ($p_\Delta$, see Theorem
\protect\ref{rejectionprob}) for the Ornstein--Uhlenbeck bridges using our
approximate method for various starting
points, $a$, and end points, $b$. The second column gives the
probability that an exact $(a,b)$-bridge is not hit by an
independent diffusion with initial distribution $p_1(b, \cdot)$. The
last column gives the probability of finding an Ornstein--Uhlenbeck a
distance $|a|$ or more from~zero}\label{OUtab1}%
\begin{tabular*}{\textwidth}{@{\extracolsep{\fill}}ld{1.2}d{1.2}d{1.3}}\hline
\multicolumn{1}{l}{$a \mapsto b$} & \multicolumn{1}{l}{Rejection prob.}
& \multicolumn{1}{l}{$ 1 - \pi$} & \multicolumn{1}{l}{$2 P(X_t > |a|)$} \\ \hline
$ -2 \mapsto-2$ & 0.09 & 0.55 & 0.05 \\
$-3 \mapsto-2$ & 0.24 & 0.74 & 0.003 \\ \hline
\end{tabular*}
\end{table}

The first row of Figure~\ref{simulation4} illustrates Theorem~\ref{diffusionbridge}. For the $(-2,-2)$ and $(-3,-2)$
Ornstein--Uhlenbeck bridges Q--Q plots compare empirical distributions
at time 0.5 of 25,000 approximate simulations to the similar empirical
distributions based on 25,000 exactly simulated Ornstein--Uhlenbeck
bridges, where the simulated bridges were removed from the sample if
the bridge was not hit by an independent diffusion with initial
distribution $p_1(b, \cdot)$. As expected from Theorem~\ref{diffusionbridge}, the two distributions appear to be equal (there is
a numerical problem in the left tail for the $(-2,-2)$-bridge). We
know from Figure~\ref{simulation3} that this distribution differs from
that of the unconditional bridge.
In the second row of Figure~\ref{simulation4}, Q--Q plots compare the
empirical distributions at time 0.5 based on 25,000 $(-2,-2)$ (left) and
$(-3,-2)$ (right) Ornstein--Uhlenbeck bridges simulated by our MH-algorithm
to the empirical distribution based on 25,000 exactly simulated
Ornstein--Uhlenbeck bridges. For the MH-algorithm the burn-in was 5000
iterations and $N=10$ (the number of $T$-values simulated in each
step). The
algorithm worked well for a much shorter burn-in and for $N=1$.
Also these distributions appear to be equal (again there is a
numerical problem in the left tail for the $(-2,-2)$-bridge). The
autocorrelations of the MH-algorithm decreases very quickly
to zero. The empirical autocorrelation function of the successive values
of the $(-3,-2)$ Ornstein--Uhlenbeck bridge at time 0.5 obtained by
our Metropolis--Hastings algorithm is plotted in Figure~\ref{simulation5}. After less than 10 iterations the correlation is
essentially zero.

Table~\ref{OUtab1} gives estimated rejection probabilities (of the
approximate rejection sampler) and the probability that an exact
$(a,b)$-bridge is not hit by an independent diffusion with initial
distribution $p_1(b, \cdot)$. As expected from the discussion above,
the probabilities of not being hit is quite
substantial. The last column gives the probability that
a stationary Ornstein--Uhlenbeck process is a distance $|a|$ or more
from zero. We see that the process will only spend very little time in
the parts of the state space, where the approximate method is biased.

%s3.2 #&#
\subsection{The hyperbolic bridge}

Next we consider the hyperbolic diffusion which is the solution to
\[
\mathrm{d}X_t= -\frac{ \theta X_t}{\sqrt{1+X_t^2}}\,\mathrm{d}t + \sigma \,\mathrm{d}W_t,
\]
with $\theta> 0$ and $\sigma>0$.
The hyperbolic diffusion was introduced by Barndorff-Nielsen \cite{oebn78}. It is
ergodic with the standardized symmetric hyperbolic distribution as
invariant distribution, see, for example, Bibby and S{\o}rensen \cite{bib-mshyp}.
In this case the transition density is not explicitly known, but
we can compare our method to the exact EA1 algorithm by Beskos, Papaspiliopoulos and Roberts \cite{Beskos2007}.
It is applicable to diffusion processes on the form
%
%e3.1 #&#
\begin{equation}
\mathrm{d}X_t = \alpha(X_t)\,\mathrm{d}t + \,\mathrm{d}W_t,
\label{eqeng}
\end{equation}
provided that $\alpha$ is continuously differentiable, and the function
$\alpha(x)^2 + \alpha'(x)$ is bounded from above and
below for all $x$, conditions satisfied by the hyperbolic diffusion
process. The algorithm by Beskos, Papaspiliopoulos and Roberts \cite{Beskos2007}
is very quick for short intervals as it essentially only requires one
simulation of a Brownian bridge if it is not rejected. Rejection in the EA1
algorithm is not very costly computationally in our example since it
is only a few points that are thrown away per rejection. Thus, we
compare our algorithm to a very efficient method.

Again we simulated 25,000 bridges using the Euler scheme with a 100
points subdivision of $[0,1]$. Also for the hyperbolic diffusion the
Euler scheme equals the Milstein scheme. The parameter values were
$\theta=
\sigma= 1$. We start with a bridge from 0 to 0 and compare our
approximate method to the exact EA1 algorithm. To the right in Figure~\ref{simulation1}, we have plotted the quantiles of the empirical
distribution at the time point $0.5$ obtained by our approximate
method against the quantiles of the empirical distribution obtained
by the exact EA1 algorithm. Also for this example the two distributions
appear to be equal.
Table~\ref{hyperbolictab} shows CPU execution times to simulate 10,000
hyperbolic diffusion bridges by our approximate method for various starting
points, $a$, and end points, $b$. Also estimated rejection probabilities
are given. The pattern is similar to that for the Ornstein--Uhlenbeck
process. For moves that are likely to
appear in data sets, the CPU times and rejection probabilities are
small, and for unlikely moves the execution time is not a problem in
applications, even though the rejection probability is quite
large. The execution time for the EA1 algorithm was 0.3 CPU seconds,
which, as expected, is faster than our method. Note that there is no
reason to consider diffusions for which the EA1 algorithm does not
work in order to compare our method to the more complicated simulation
methods EA2 and EA3 in Beskos, Papaspiliopoulos and Roberts \cite{Beskos2007} and Beskos \textit{et al}. \cite{beskos7}. The EA2
and EA3 algorithms are clearly more time consuming than EA1, while
execution times for our methods can be expected to be approximately
as for the two examples considered here.
%
%t3 #&#
\begin{table}
\tablewidth=\textwidth
\tabcolsep=0pt
\caption{The CPU execution time (in seconds) used to simulate 10,000
hyperbolic diffusion bridges with $\theta= \sigma= 1$ by our
method for various starting points, $a$, and end points, $b$. Also
estimated rejection probabilities are given ($p_\Delta$, see Theorem
\protect\ref{rejectionprob})}
\label{hyperbolictab}%
\begin{tabular*}{\textwidth}{@{\extracolsep{\fill}}ld{2.1}d{1.2}}\hline
\multicolumn{1}{l}{$a \mapsto b$} & \multicolumn{1}{l}{CPU (sec.)} & \multicolumn{1}{l}{Rejection prob.} \\ \hline
$\hphantom{-}{0 \mapsto0}$ & 0.6 & 0.14 \\
$\hphantom{-}{0 \mapsto1}$ & 0.8 & 0.36 \\
$\hphantom{-}{0 \mapsto2}$ & 2.1 & 0.77 \\
$-1 \mapsto1$ & 2.0 & 0.76 \\
$-1 \mapsto2$ & 12.6 & 0.96 \\ \hline
\end{tabular*}
\end{table}
%
%t4 #&#
\begin{table}[b]
\tablewidth=\textwidth
\tabcolsep=0pt
\caption{The CPU execution time (in seconds) used to simulate 10,000
hyperbolic $(0,0,\Delta,0)$-bridges with $\theta= \sigma^2 = 4$
for our approximate method
and for the EA1 method in Beskos \textit{et al}. \cite{beskos7} for different
interval lengths $\Delta$. Also the number of rejections while
simulating the 10,000 trajectories is given}
\label{intervallenght}
\begin{tabular*}{\textwidth}{@{\extracolsep{\fill}}d{1.1}d{1.2}d{3.0}d{3.2}d{9.0}}\hline
& \multicolumn{2}{l}{Present paper} &
\multicolumn{2}{l}{Beskos \textit{et al}. \cite{beskos7} EA1} \\ [-5pt]
& \multicolumn{2}{c}{\hrulefill} &
\multicolumn{2}{c}{\hrulefill} \\
\multicolumn{1}{l}{$\Delta$} & \multicolumn{1}{l}{CPU time} & \multicolumn{1}{l}{$\#$ rejections} &
\multicolumn{1}{l}{CPU time}
& \multicolumn{1}{l}{$\#$ rejections} \\ \hline
0.5 & 0.52 & 819 & 0.28 & 14\mbox{,}497 \\
1.0 & 0.99 & 307 & 0.59 & 53\mbox{,}087 \\
1.5 & 1.45 & 102 & 1.05 & 163\mbox{,}599 \\
2.0 & 1.93 & 44 & 1.92 & 457\mbox{,}226 \\
2.5 & 2.40 & 17 & 4.00 & 1\mbox{,}242\mbox{,}922 \\
3.0 & 2.88 & 6 & 10.01 & 3\mbox{,}491\mbox{,}838 \\
3.5 & 3.36 & 2 & 26.86 & 9\mbox{,}357\mbox{,}310 \\
4.0 & 3.83 & 0 & 75.79 & 25\mbox{,}232\mbox{,}418 \\
4.5 & 4.31 & 0 & 222.09 & 69\mbox{,}299\mbox{,}642 \\
5.0 & 4.79 & 0 & 641.70 & 187\mbox{,}069\mbox{,}771 \\ \hline
\end{tabular*}
\end{table}

\begin{figure}

\includegraphics{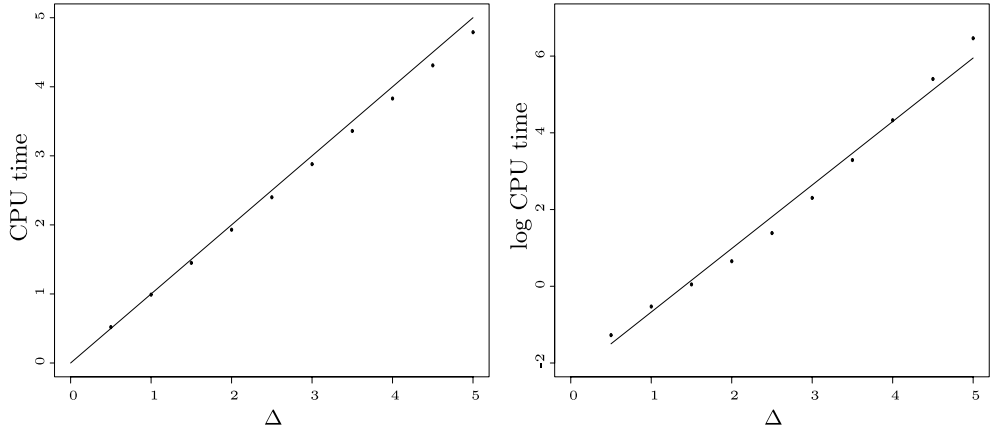}

\caption{The CPU execution time (in seconds) used to simulate 10,000
hyperbolic $(0,0,\Delta,0)$-bridges with $\theta= \sigma^2 = 4$.
In the left plot the CPU time is
plotted against $\Delta$ for our approximate method, while in the
right plot the logarithm of the CPU time is plotted against $\Delta$
for the EA1 method in Beskos \textit{et al}. \cite{beskos7}.}\label{intervallengthplots}
\end{figure}

Beskos, Papaspiliopoulos and Roberts \cite{Beskos2007} noted that the computing time of their exact
algorithm is large for diffusion bridges over long time intervals. It is
therefore of interest to compare computer time and rejection
probabilities for our approximate algorithm to the EA1 algorithm.
To do so, we simulated 10,000 trajectories of the $(0,0,\Delta,0)$-bridge
for the hyperbolic diffusion with $\theta= \sigma^2 = 4$. This was
done for values of the interval length $\Delta$ ranging from 0.5 to
5. The CPU execution time (in seconds) used to simulate the 10,000
trajectories are given in Table~\ref{intervallenght}. We see that for
this particular diffusion the two methods use the same CPU time for an
interval length of two. For smaller interval lengths the exact algorithm
is somewhat faster, whereas our approximate method
is much faster for long intervals. The simulations confirm that the
computational complexity of the proposed method is linear in the
interval length $\Delta$ as shown in Section~\ref{sec2}, whereas the complexity
appears to
grow at least exponentially with $\Delta$ for the exact algorithm; see
Figure~\ref{intervallengthplots}. The main reason is that for long
intervals the number of rejections becomes very large for the
algorithm in Beskos, Papaspiliopoulos and Roberts \cite{Beskos2007}, while our approximate algorithm
has a very small rejection probability for long intervals. The rapid
decrease of the rejection probabilities for the approximate method as a
function of $\Delta$ is expected from Theorem~\ref{rejectionprob}.
%
%f7 #&#

%s4 #&#
\section{Maximum likelihood estimation}\label{sec4}

The main motivation for the theory developed in this papers is the
central role diffusion bridge simulation plays in simulated
likelihood-based inference for processes of the diffusion
type. Therefore, we end the paper by giving two examples of application
of our diffusion bridge simulation method to maximum likelihood
estimation first for discretely observed diffusion processes then for
integrated diffusions observed with measurement error.

%s4.1 #&#
\subsection{Discretely observed diffusions}\label{sec4.1}

First, we present an EM-algorithm for finding the maximum likelihood
estimator for discretely observed diffusion processes. We also briefly
discuss aspects of Bayesian inference.

Consider the diffusion process
%
%e4.1 #&#
\begin{equation}
\label{basicmodel} \mathrm{d}X_t = b_{\alpha}(X_t)\,\mathrm{d}t+
\sigma_\beta(X_t)\,\mathrm{d}W_t ,
\end{equation}
where $\alpha$ and $\beta$ are unknown parameters to be estimated, and
$W$ is the standard Wiener process. We assume that $\sigma_\beta
(x)>0$ for all $x$ in the state interval. Suppose that the only data
available from a realization of the diffusion process are observations
at times $t_1<t_2<\cdots <t_n$, $x_i = X_{t_i}$, $i=1,\ldots ,n$.

As explained in the \hyperref[intro]{Introduction}, discrete time observation of a
continuous time process can be viewed as an incomplete observation
problem, so the EM-algorithm (Dempster, Laird and Rubin \cite{dempster}) is a natural method
for finding the maximum likelihood estimator of the parameters. Maximum
likelihood estimation for discretely observed Markov jump processes
was treated in this way by Bladt and S{\o}rensen \cite{mbms,mbms09}.
Unfortunately, the probability measures corresponding to
complete continuous time observation of the diffusion model given by
(\ref{basicmodel}) are singular because the diffusion coefficient
depends on the parameter $\beta$. It is therefore not straightforward
to implement the EM-algorithm, but an approach in the spirit of
Roberts and Stramer \cite{roberts} was proposed by Beskos, Papaspiliopoulos and Roberts \cite{Beskos2006}. In the following,
we summarize a modification of this approach using our diffusion
bridge simulation technique.

The transformation
%
%e4.2 #&#
\begin{equation}
\label{h-trans} h_{\beta}(x)=\int_{x^{\ast}}^x
\frac{1}{\sigma_\beta(y)}\,\mathrm{d}y
\end{equation}
is essential. Here $x^{\ast}$ is some arbitrary, but appropriately
chosen, point in the state interval. By Ito's formula,
$Y_t=h_\beta(X_t)$ solves
%
%e4.3 #&#
\begin{equation}
\label{Y} \mathrm{d}Y_t=\mu_{\alpha,\beta}(Y_t)\,\mathrm{d}t+\mathrm{d}W_t,
\end{equation}
where
\[
\mu_{\alpha,\beta}(y)=\frac{b_{\alpha}(h_{\beta}^{-1}(y))} {
\sigma_{\beta}(h_\beta^{-1}(y))}-\frac{1}{2}\sigma_\beta^{\prime}
\bigl(h_\beta^{-1}(y) \bigr).
\]
In (\ref{Y}), the diffusion coefficient does not depend on the
parameters, so the probability measures are equivalent and the
likelihood function can be found. To do so the function
%
%e4.4 #&#
\begin{equation}
\label{g} g_{\alpha,\beta}(x)= s_{\alpha,\beta}(x) \tfrac12 \log\bigl(
\sigma_\beta(x)\bigr),
\end{equation}
where
%
%e4.5 #&#
\begin{equation}
\label{s} s_{\alpha,\beta}(x) = \int_{x^{\ast}}^x
\frac{b_\alpha(z)}{\sigma_\beta^2(z)}\,\mathrm{d}z,
\end{equation}
is needed. Note that $\int_{y^{\ast}}^y \mu_{\alpha,\beta}(z)\,\mathrm{d}z =
g_{\alpha,\beta}(h^{-1}_\beta(y))-g_{\alpha,\beta}(h^{-1}_\beta
(y^{\ast}))$,
and that the functions $g_{\alpha,\beta}$ and
$s_{\alpha,\beta}$ are closely related to the density
$\varphi_{\alpha,\beta}$ of the stationary distribution of the
original diffusion model given by (\ref{basicmodel}). Specifically,
$s_{\alpha,\beta}(x)$ equals $\frac12 \log(\sigma_\beta(x)^2
\varphi_{\alpha,\beta} (x))$ apart from an additive constant.
Thus, when the stationary density is known, the only problem is to find
$h_\beta$ and its inverse. This is for instance the case for the
Pearson diffusions studied by Forman and S{\o}rensen \cite{julie}.

The problem with the transformation $h_{\beta}$ is that it is
parameter dependent, while we need to keep the original discrete time
data fixed when running the EM-algorithm. To get around this problem,
define
\[
Y_t^{*}(\beta,\beta_0)=Z^{(i,\alpha_0,\beta_0)}_t
+\frac{(t_i-t)(h_\beta(x_{i-1})-h_{\beta
_0}(x_{i-1}))+(t-t_{i-1})(h_\beta
(x_i)-h_{\beta_0}(x_i))}{t_i-t_{i-1}}
\]
for $t_{i-1} \leq t \leq t_i$, $i=2, \ldots, n$. Here
$Z^{(i,\alpha_0,\beta_0)}_t$ denotes the
$(t_{i-1},h_{\beta_0}(x_{i-1}),t_i,h_{\beta_0}(x_i))$-bridge for the
diffusion (\ref{Y}) with parameter values $\alpha_0$ and $\beta_0$,
and $Z^{(i,\alpha_0,\beta_0)}_t$, $i=2, \ldots, n$ are independent.
Then the EM-algorithm works as follows. Let $\alpha_0,\beta_0$ be
initial values of the parameters.
\begin{longlist}[(3)]
\item[(1)] (E-step) Calculate the function
\begin{eqnarray*}
q(\alpha,\beta) &=& g_{\alpha,\beta}(x_n)-g_{\alpha,\beta}(x_1)\\
&&{}-\frac12 \sum_{i=2}^n
\bigl[h_\beta(x_i)-h_\beta (x_{i-1})
\bigr]^2/(t_i-t_{i-1})
- \sum_{i=2}^n \log\bigl(
\sigma_\beta(x_i)\bigr) \\
&&{}-\frac{1}{2} \sum
_{i=2}^n \Exp_{Z^{(i,\alpha_0,\beta_0)}} \biggl( \int
_{t_{i-1}}^{t_i} \bigl[ \mu_{\alpha,\beta}^{\prime}
\bigl(Y^{*}_t(\beta,\beta_0)\bigr)+
\mu_{\alpha
,\beta} \bigl(Y^{*}_t(\beta,
\beta_0)\bigr)^2 \bigr] \,\mathrm{d}t \biggr).
\end{eqnarray*}
\item[(2)] (M-step) $(\alpha_0,\beta_0)=\argmax_{\alpha
,\beta}
q(\alpha,\beta)$.
\item[(3)] GO TO (1).
\end{longlist}
In the E-step, $\Exp_{Z^{(i,\alpha_0,\beta_0)}}$ means that the data
points are fixed so that only the diffusion bridge is random, and
expectation is with respect to the distribution of the diffusion
bridge. Thus, the expectations in the E-step can be approximated
by simulating diffusion bridges by our exact MH-method and averaging
(after a burn-in period).
Arguments that $q(\alpha,\beta)$ is the conditional expectation of the
relevant continuous time likelihood function can be found in
Roberts and Stramer \cite{roberts} and Beskos, Papaspiliopoulos, Roberts and Fearnhead \cite{Beskos2006}. As pointed out in the latter
paper, the conditional expectation can also be calculated as
\[
\Exp_{Z^{(i,\alpha_0,\beta_0)},U} \bigl( \mu_{\alpha,\beta
}^{\prime}
\bigl(Y^{*}_U(\beta,\beta_0)\bigr)+
\mu_{\alpha,\beta} \bigl(Y^{*}_U(\beta,
\beta_0)\bigr)^2 \bigr),
\]
where $U$ is a uniformly distributed random variable on
$[t_{i-1},t_i]$ that is independent of $Z^{(i,\alpha_0,\beta_0)}_t$
and the data.

In the M-step the maximization of $q(\alpha,\beta)$ must in general
be done by a suitable maxi\-mization algorithm. With modern software
(e.g., the R-function optim), this is not a problem. When the drift
of the original diffusion model (\ref{basicmodel}) depends linearly
on the vector of parameters $\alpha$, that is, when
%
%e4.6 #&#
\begin{equation}
\label{paramlindrift} b_\alpha(x) = \alpha_1
a_1(x) + \cdots+ \alpha_k a_k(x),
\end{equation}
where $a_1, \ldots, a_k$ are known functions, then
the maximization problem is simplified somewhat. When the drift has this
form, and when the diffusion parameter $\beta$ is fixed, the
model for continuous time observation of $X$ as well as the
transformed process $Y$ is an exponential family of stochastic
processes, see K\"{u}chler and S{\o}rensen \cite{kuso}, page 27. We can therefore take advantage of well
known properties of exponential families of diffusions.

For the EM-algorithm the specification (\ref{paramlindrift}) implies
that the function $q(\alpha,\beta)$ has the form
\[
q(\alpha,\beta) = \sum_{i=1}^k
\alpha_i H_{i,\beta} - \frac12 \sum
_{i=1}^k \sum_{j=1}^k
\alpha_i \alpha_j B_{i,j,\beta} + G_\beta,
\]
where
\begin{eqnarray*}
H_{i,\beta} &=& s_{i,\beta}(x_n) - s_{i,\beta}(x_1)
\\
&& {} + \sum_{j=2}^n
\Exp_{Z^{(j,\alpha_0,\beta_0)}} \biggl( \int_{t_{j-1}}^{t_j} \biggl[
a_i\bigl(h_\beta^{-1}\bigl(Y^*_t(
\beta,\beta_0)\bigr)\bigr) (\log\sigma_\beta)'
\bigl(h_\beta^{-1}\bigl(Y^*_t(\beta,
\beta_0)\bigr)\bigr) \\
&&\hphantom{{}+ \sum_{j=2}^n
\Exp_{Z^{(j,\alpha_0,\beta_0)}} ( \int_{t_{j-1}}^{t_j} \biggl[}{}- \frac12 a_i'
\bigl(h_\beta^{-1}\bigl(Y^*_t(\beta,
\beta_0)\bigr)\bigr) \biggr] \,\mathrm{d}t \biggr),
\end{eqnarray*}
with $s_{i,\beta}(x) = \int_{x^\ast}^x a_i(y)/\sigma^2_\beta(y)\, \mathrm{d}y$,
\[
B_{i,j,\beta} = \sum_{j=2}^n
\Exp_{Z^{(j,\alpha_0,\beta_0)}} \biggl( \int_{t_{j-1}}^{t_j}
\frac{a_i(h_\beta^{-1}(Y^*_t(\beta,\beta_0)))a_j(h_\beta^{-1}
(Y^*_t(\beta,\beta_0)))}{\sigma^2_\beta(h_\beta^{-1}
(Y^*_t(\beta,\beta_0)))} \,\mathrm{d}t \biggr),
\]
and
\begin{eqnarray*}
G_\beta&=& -\frac12 \log \bigl( \sigma_\beta(x_n)/
\sigma_\beta(x_1) \bigr)\\
&&{} -\frac12 \sum
_{i=2}^n \bigl[h_\beta(x_i)-h_\beta(x_{i-1})
\bigr]^2/(t_i-t_{i-1}) - \sum
_{i=2}^n \log\bigl(\sigma_\beta(x_i)
\bigr)
\\
&& {} + \frac14 \sum_{j=2}^n
\Exp_{Z^{(j,\alpha_0,\beta_0)}} \biggl( \int_{t_{j-1}}^{t_j} \biggl[
\sigma_\beta''\bigl(h_\beta^{-1}
\bigl(Y^*_t(\beta,\beta_0)\bigr)\bigr)
\sigma_\beta \bigl(h_\beta^{-1}
\bigl(Y^*_t(\beta,\beta_0)\bigr)\bigr) \\
&&\hphantom{{}+ \frac14 \sum_{j=2}^n
\Exp_{Z^{(j,\alpha_0,\beta_0)}} ( \int_{t_{j-1}}^{t_j} \biggl[}{}- \frac12 \bigl\{
\sigma_\beta'\bigl(h_\beta^{-1}
\bigl(Y^*_t(\beta,\beta_0)\bigr)\bigr)\bigr
\}^2 \biggr] \,\mathrm{d}t \biggr).
\end{eqnarray*}
For a fixed value of $\beta$, the function
$\alpha\mapsto q(\alpha,\beta)$ is maximal for
\[
\hat{\alpha}(\beta) = \bfB_\beta^{-1} \bfH_\beta,
\]
where $\hat{\alpha} = (\hat{\alpha}_1, \ldots, \hat{\alpha}_k)^T$,
$\bfH_\beta= (H_{1,\beta}, \ldots, H_{k,\beta})^T$ and $\bfB_\beta= \{
B_{i,j,\beta} \}$. This is provided that $\bfB_\beta$ is invertible,
which it is when the functions $a_i$, $i=1,\ldots,k$ are linearly
independent. Thus, $q(\alpha,\beta)$ attains its maximal value at
$(\hat{\alpha}(\hat{\beta}),\hat{\beta})$, where $\hat{\beta}$ maximizes
\[
\beta\mapsto q\bigl(\hat{\alpha}(\beta),\beta\bigr) = \tfrac12
H_\beta^T B_\beta^{-1}
H_\beta+ G_\beta.
\]

The Gibbs sampler for Bayesian inference for discretely observed
diffusion processes proposed by Roberts and Stramer \cite{roberts} can also be modified by
replacing the MCMC algorithm for simulating diffusion bridges in that
paper by our diffusion bridge simulation method. We will not go into
any detail for general diffusions, but will limit ourselves to
pointing out that when the drift has the form (\ref{paramlindrift}), then
the (continuous time) posterior distribution of $\alpha$ simplifies.
Choose as the prior for $\alpha$ the conjugate prior for an exponential
family of diffusions (see K\"{u}chler and S{\o}rensen \cite{kuso}, page 51), which here is a multivariate
normal distribution
with expectation $\bar{\alpha}$ and covariance matrix $\bfSigma$. Then
the posterior of $\alpha$ (given $\beta=\beta_0$ and given simulated
diffusion
bridges) is a $k$-dimensional normal distribution with
expectation $(\bfSigma^{-1} + \tilde{\bfB}_{\beta_0})^{-1}
(\bfSigma^{-1}\bar{\alpha} + \tilde{\bfH}_{\beta_0})$ and
covariance matrix
$(\bfSigma^{-1} + \tilde{\bfB}_{\beta_0})^{-1}$, where $\tilde
{\bfH}_\beta
= (\tilde{H}_{1,\beta}, \ldots, \tilde{H}_{k,\beta})^T$, $\bfB_\beta= \{
B_{i,j,\beta} \}$,
\begin{eqnarray*}
\tilde{H}_{i,\beta} &=& s_{i,\beta}(x_n) -
s_{i,\beta}(x_1)
\\
&& {} + \sum_{i=2}^n \int
_{t_{i-1}}^{t_i} \biggl[ a_i
\bigl(h_\beta^{-1} \bigl(Y^*_t(\beta,
\beta_0)\bigr)\bigr) (\log\sigma_\beta)'
\bigl(h_\beta^{-1}\bigl(Y^*_t(\beta,
\beta_0)\bigr)\bigr)\\
&&\hphantom{{} + \sum_{i=2}^n \int
_{t_{i-1}}^{t_i} \biggl[}{} - \frac12 a_i'
\bigl(h_\beta^{-1}\bigl(Y^*_t(\beta,
\beta_0)\bigr)\bigr) \biggr] \,\mathrm{d}t,
\end{eqnarray*}
and
\[
\tilde{B}_{i,j,\beta} = \sum_{i=2}^n
\int_{t_{i-1}}^{t_i} \frac{a_i(h_\beta^{-1}(Y^*_t(\beta,\beta_0)))
a_j(h_\beta^{-1}(Y^*_t(\beta,\beta_0)))} {
\sigma^2_\beta(h_\beta^{-1}(Y^*_t(\beta,\beta_0)))} \,\mathrm{d}t.
\]

%s4.2 #&#
\subsection{Integrated diffusions observed with measurement error}

Here we present an EM-algorithm to find the maximum likelihood
estimator when an integrated diffusion is observed with measurement
errors. The method was proposed and studied by Baltazar-Larios and S{\o}rensen \cite{fernando}.

We consider again the diffusion process $X$ given by (\ref{basicmodel}),
but here the data are of the form
%
%e4.7 #&#
\begin{equation}
\label{2.1} V_i=\int_{t_{i-1}}^{t_i}X_{s}\,\mathrm{d}s+Z_i,\qquad
i=1,\ldots,n,
\end{equation}
where $Z_i\sim N(0,\tau^2)$, $i=1,\ldots,n$ are mutually independent and
independent of ${X}$, and $t_0 = 0$. We assume that $X$ is stationary
and ergodic. The variance of the
measurement error, $\tau^2$, is an extra unknown parameter, so we
need to estimate the parameter $\theta=(\alpha, \beta,\tau^2)$.
We can think of the data set $V=(V_1,\ldots,V_n)$ as an incomplete
observation of a full data set given by the sample path $X_t, t \in
[0,t_n]$ and the measurement errors $Z_1, \ldots, Z_n$, or
equivalently $X_t, t \in[0,t_n]$ and $V=(V_1,\ldots,V_n)$. To apply
the EM-algorithm, we need to find the likelihood function for the full
data set and the conditional expectation of this full log-likelihood
function given the observations $V=(V_1,\ldots,V_n)$.

Conditionally on the sample path of $X$, the observations $V_i$,
$i=1,\ldots,n$ are independent and normal distributed with expectation
$\int_{t_{i-1}}^{t_i}X_s\,\mathrm{d}s$ and variance $\tau^2$. Again we need to
apply the\vspace*{1pt} transformation (\ref{h-trans}) because the probability
measures are singular. By expressing the data $V_i$ in terms of the
process $Y$, using the parameter-dependent transformation $h_\beta$,
and by Girsanov's theorem, we find that the log-likelihood function
for the full data set is
%
%e4.8 #&#
\begin{eqnarray}
\label{fullloglik} &&\ell\bigl(\theta;V_1,\ldots,V_n,Y_t,t
\in[0,t_n]\bigr)\nonumber\\
 &&\quad = \sum_{i=1}^n
\log \varphi\biggl(V_i;\int_{t_{i-1}}^{t_i}h^{-1}_\beta(Y_s)\,\mathrm{d}s,
\tau^2\biggr)+ g_{\alpha,\beta}\bigl(h^{-1}_\beta(Y_{t_n})
\bigr)
\\
&& \qquad {}  -g_{\alpha,\beta}\bigl(h^{-1}_\beta(Y_{0})
\bigr) -\frac{1}{2}\int_{0}^{t_n} \bigl(
\mu_{\alpha,\beta}(Y_t)^2 +\mu'_{\alpha,\beta}(U_t)
\bigr)\,\mathrm{d}t,
\nonumber
\end{eqnarray}
where $\varphi(x;a,b)$ denotes the normal density function with mean
$a$ and variance $b$. The EM-algorithm works as follow. Let the
initial value $\hat{\theta} = (\hat\alpha, \hat\beta, \hat\tau^2)$
be any value of the parameter vector $\theta$.
\begin{longlist}
\item[(1)] (E-step) Generate $M$ sample paths of the diffusion
process $X$,
$X^{(k)},k=1,\ldots,M$, conditional on the observations
$V_1,\ldots,V_n$ using the parameter value $\hat{\theta}$, and calculate
\[
g(\theta)=\frac{1}{M-M_0}\sum_{k=M_0+1}^M
\ell\bigl(\theta;Y_1,\ldots,Y_n,h_{\hat\beta}
\bigl(X_t^{(k)}\bigr),t\in[0,t_n]\bigr)
\]
for a suitable burn-in period $M_0$ and $M$ sufficiently large.
\item[(2)] (M-step) $\hat{\theta}=$ argmax $g(\theta).$
\item[(3)] GO TO (1).
\end{longlist}

To implement
this algorithm, the main issue is how to generate sample paths of $X$
conditionally on $V_1,\ldots,V_n$, where the relation between the
$V_i$s and $X$ is given by (\ref{2.1}). This can be done by means of a
Metropolis--Hastings algorithm. However, if the sample path in the
entire time interval $[0,t_n]$ is updated in one step, the rejection
probability is typically very large. Therefore, it is more efficient to
randomly divide the time interval into subintervals and update the
sample path in each of the subintervals conditional on the rest of the
sample path, which corresponds to simulating a diffusion
bridge in each subinterval (except the end-intervals). This is a
modification of the method in Chib \textit{et al}. \cite{chibpittshephard}, where we use
the algorithm for diffusion bridge simulation proposed in Section~\ref{sec2}.
In the following, the parameter value $\theta$ is fixed. Start by
generating an initial unrestricted stationary sample path,
$\{X_t^{(0)}\dvt  t\in[0,t_n]\}$, of the diffusion given by
(\ref{basicmodel}), and set $l=1$.
\begin{longlist}[(3)]
\item[(1)] Generate a sample path $\{X_t^{(l)}\dvt  t\in[0,t_n]\}$ conditional
on $Y$ by updating subsets of the sample path:
\begin{enumerate}[(b)]
\item[(a)] Randomly split the time interval from $[0,t_n]$ into $K$
blocks, and write these subsampling times as
$0=\tau_0\leq\tau_1\leq\cdots\leq\tau_K=t_n,$
where each $\tau_i$ is one of the end-points of the integration
intervals, $t_j$, $j=0,\ldots,n$. Let $Y_{\{k\}}$ denote the collection
of all observations $Y_j$ for which $\tau_{k-1} < t_j \leq\tau_k$.
\item[(b)] Draw $X_0^{(l)}$ from the stationary distribution,
and simulate the conditional subpath $\{X_t^{(l)}\dvt
t\in[\tau_{k-1},\tau_k]\} | Y_{\{k\}},
X_{\tau_{k-1}}^{(l)},X_{\tau_k}^{(l-1)},$
for $k=1,\ldots,K-1$. Finally, simulate
$\{X_t^{(l)}\dvt  t\in[\tau_{K-1},\allowbreak \tau_K]\} | Y_{\{K\}},
X_{\tau_{K-1}}^{(l)}$.
\end{enumerate}
\item[(2)] $l=l+1$.
\item[(3)] GO TO (1).
\end{longlist}

The random time intervals can for instance be generated by independent
Poisson variables. Simulation of a $(\tau_{k-1},a,\tau_k,b)$-bridge
conditional on $Y_{\{k \}}$, the data in $(\tau_{k-1}, \tau_{k}]$,
can be done by a Metropolis--Hastings algorithm that uses the bridge
simulation method introduced in Section~\ref{sec2} as proposal and accepts a
proposed bridge $X^{(l)}$ with probability
\[
\min \Biggl( 1,\prod_{i=1}^{n_k}
\frac{\varphi(Y_{j+i}; \int_{t_{j+i-1}}^{t_{j+i}} X_s^{(l)}\,\mathrm{d}s,\tau^2)} {
\varphi(Y_{j+i}; \int_{t_{j+i-1}}^{t_{j+i}}X_s^{(l-1)}\,\mathrm{d}s,\tau^2)} \Biggr),
\]
where the end-point $\tau_{k-1}$ is equal to $t_j$, $n_k$ is the
number of observations in the interval $(\tau_{k-1}, \tau_{k}]$
(the observations are $Y_{j+1}, \ldots, Y_{j+n_k}$), and $X^{(l-1)}$
is the bridge in the previous step of the MH-algorithm.

In a simulation study by Baltazar-Larios and S{\o}rensen \cite{fernando}, where the approximate
algorithm of Section \ref{sec2.1} was used to simulate diffusion bridges,
this EM-algorithm worked well. The study considered 1500 observations of
Ornstein--Uhlenbeck and CIR processes integrated over time intervals of
length one. Three versions of the EM-algorithm were investigated,
where the expected number of observations in each random subinterval
was 11, 21 and 31, respectively.

%s5 #&#
\section{Conclusion}\label{sec5}

We have presented a straightforward way of simulating an approximation
to a diffusion bridge and an easily implementable Metropolis--Hastings
algorithm that uses the approximate simulation as proposal and has
exact diffusion bridges as the target distribution. Advantages of the
new method is that it is easy to understand and to implement, that the
same simple algorithm can be used for all one-dimensional diffusion
processes with finite speed measure, and most importantly, that the
method works particularly well for long time intervals, where other
methods tends not to work or to be very time consuming. The method
allows the use of simple simulation procedures like the Milstein
scheme for bridge simulation. The simulation study showed that the
one-dimensional distributions obtained by the approximate method
compare accurately to the results from exact simulations for bridges
corresponding to data that are likely in discrete-time samples from
diffusion models.

For ergodic diffusions, the computational complexity was shown to be
linear in the interval length for our approximate method as well as
for our exact method. The simulation study showed that the computing
time for the approximate algorithm for small time intervals is of the
same order of magnitude as for the exact EA1 method and for long
intervals it is much faster than the EA1 method. Thus our new diffusion
bridge simulation method is highly suitable for likelihood inference
for discretely observed diffusions and can be used to simplify and in
some cases speed up methods for likelihood inference and Bayesian
inference (the EM-algorithm and the Gibbs sampler) for discretely
observed diffusion processes. The method is also potentially useful
for inference for more general diffusion type models like stochastic
volatility models.


\section*{Acknowledgements}

The authors are grateful to the referees for
their insightful comments that have improved the paper very considerably.
The research of Michael S\o rensen was supported by the Danish
Center for Accounting and Finance funded by the Danish Social
Science Research Council, by the Center for Research in
Econometric Analysis of Time Series funded by the Danish National
Research Foundation, and by a grant from the University of Copenhagen
Programme of Excellence.

%suskaldyti doi

% imsref loaded by jurgita.kaciuliene, 2013-02-04 15:10:10

\printhistory

\end{document}